\documentclass[11pt]{article}

\usepackage{amssymb}
\usepackage{amsmath,amsfonts}
\usepackage{amsthm}
\usepackage{latexsym}
\usepackage{mathrsfs}
\usepackage{color}

\usepackage[a4paper]{geometry}
\newcommand{\beq}{\begin{equation}}
\newcommand{\eeq}{\end{equation}}
\newcommand{\beqa}{\begin{eqnarray}}
\newcommand{\eeqa}{\end{eqnarray}}
\newcommand{\beqas}{\begin{eqnarray*}}
\newcommand{\eeqas}{\end{eqnarray*}}
\newcommand{\ba}{\begin{array}}
\newcommand{\ea}{\end{array}}
\newcommand{\bi}{\begin{itemize}}
\newcommand{\ei}{\end{itemize}}
\newcommand{\nn}{\nonumber}

\DeclareMathOperator*{\Argmin}{Argmin}  

\newcommand{\mcK}{{\mathcal K}}

\newcommand{\mcL}{{\mathcal L}}

\newcommand{\mcN}{{\mathcal N}}
\newcommand{\mcJ}{{\mathcal J}}
\newcommand{\mcT}{{\mathcal T}}
\newcommand{\mcI}{{\mathcal I}}
\newcommand{\mcS}{{\mathcal S}}

\newcommand{\prox}{\mathrm{prox}}
\newcommand{\dom}{\mathrm{dom}}
\newcommand{\dist}{\mathrm{dist}}
\newcommand{\argmin}{\arg\min}

\newcommand{\ceil}[1]{\left\lceil #1 \right\rceil}
\newcommand{\floor}[1]{\left\lfloor #1 \right\rfloor}

\newcommand{\Inner}[2]{\left\langle #1, #2 \right\rangle}

\newtheorem{lemma}{Lemma}
\newtheorem{thm}{Theorem}

\newtheorem{defi}{Definition}
\newtheorem{prop}{Proposition}

\newtheorem{ass}{Assumption}

\newtheorem{rem}{Remark}

\newtheorem{algo}{Algorithm}

\newcounter{spb}
\def\bB{{\bar B}}
\def\bC{{\bar C}}
\def\bD{{\bar D_{\Lambda}}}
\def\cA{{\cal A}}
\def\cO{{\cal O}}
\def\tlambda{{\tilde \lambda}}
\def\tx{{\tilde x}}

\begin{document}
\title{Iteration Complexity of First-Order Augmented Lagrangian Methods \\ for Convex Conic Programming}
\author{
	Zhaosong Lu%
	\thanks{
		Department of Industrial and Systems Engineering, University of Minnesota, USA (email: {\tt zhaosong@umn.edu}). The work of this author was partially supported by NSF Award IIS-2211491.
	} 
	\and
	Zirui Zhou%
	\thanks{Huawei Technologies Canada, Burnaby, BC, Canada (email: {\tt zirui.zhou@huawei.com}).}
}
\date{March 8, 2021 (Revised: February 11, 2022; July 2, 2022)}
\maketitle

\begin{abstract}
In this paper we consider a class of convex conic programming. In particular, we first propose an inexact augmented Lagrangian (I-AL) method that resembles the classical I-AL method for solving this problem, in which the augmented Lagrangian subproblems are solved approximately by a variant of Nesterov's optimal first-order method. We show that the total number of first-order iterations of the proposed I-AL method for finding an $\epsilon$-KKT solution is at most $\mathcal{O}(\epsilon^{-7/4})$. We then propose an adaptively regularized I-AL method and show that it achieves a first-order iteration complexity $\mathcal{O}(\epsilon^{-1}\log\epsilon^{-1})$, which significantly improves existing complexity bounds achieved by first-order I-AL methods for finding an $\epsilon$-KKT solution. Our complexity analysis of the I-AL methods is based on a sharp analysis of inexact proximal point algorithm (PPA) and the connection between the I-AL methods and inexact PPA. It is vastly different from existing complexity analyses of the first-order I-AL methods in the literature, which typically regard the I-AL methods as an inexact dual gradient method.

\end{abstract}

\bigskip


\bigskip


\section{Introduction}\label{sec:intro}
In this paper we consider convex conic programming in the form of
\begin{equation}\label{eq:conic-cp}
\begin{array}{rl}
F^*=\min & \left\{F(x) := f(x) + P(x)\right\} \\
\mbox{s.t.} &  g(x)\preceq_\mcK 0,
\end{array}
\end{equation}
where $f, P:\Re^n\rightarrow(-\infty,\infty]$ are proper closed convex functions, $\mcK$ is a closed convex cone in $\Re^m$, the symbol $\preceq_\mcK$ denotes the partial order induced by $\mcK$, that is, $y\preceq_\mcK z$ if and only if $z-y\in\mcK$, and the mapping $g:\Re^n\rightarrow\Re^m$ is convex with respect to $\mcK$, that is,
\begin{equation}
\label{eq:conic-cvx}
g(\alpha x + (1-\alpha)y) \preceq_\mcK \alpha g(x) + (1-\alpha)g(y), \quad\forall x,y\in\Re^n, \; \alpha\in[0,1].
\end{equation}
The associated Lagrangian dual problem of \eqref{eq:conic-cp} is given by 
\beq \label{dual-prob}
d^*=\sup\limits_{\lambda\in\mcK^{*}} \inf\limits_x\  \{ f(x) + P(x) + \Inner{\lambda}{g(x)}\},
\eeq
where $\mcK^*$ is the dual cone of $\mcK$ (see Section \ref{sec:preliminaries}). 
We make the following additional assumptions on problems \eqref{eq:conic-cp} and \eqref{dual-prob} throughout this paper.
\begin{ass}
	\label{ass:function}
	\begin{itemize}
		\item[(a)] The proximal operator associated with $P$ can be exactly evaluated, and the domain of $P$, denoted by $\dom(P)$, is compact.\footnote{Some problems with unbounded $\dom(P)$ can be reformulated as the ones satisfying this assumption. For example, when $f$ is bounded below and $P$ is coercive, the problem generally can be reformulated as the one with the objective $f+\tilde P$ for some $\tilde P$ with a compact domain. Such a problem often arises in sparse or low-rank learning, in which $f$ is typically a nonnegative loss function and $P$ is the $\ell_1$- or nuclear-norm.}
		\item[(b)] The projection onto $\mcK$ can be exactly evaluated.
		\item[(c)] The functions $f$ and $g$ are continuously differentiable on an open set $\Omega$ containing $\dom(P)$,  and $\nabla f$ and $\nabla g$  are Lipschitz continuous on $\Omega$ with Lipschitz constants $L_{\nabla f}$ and $L_{\nabla g}$, respectively.\footnote{The symbol $\nabla g$ denotes the transpose of the Jacobian of $g$.}
\item[(d)] The strong duality holds for problems \eqref{eq:conic-cp} and \eqref{dual-prob}, that is, both problems have optimal solutions and moreover their optimal values $F^*$ and $d^*$ are equal. 
	\end{itemize}
\end{ass}

Problem \eqref{eq:conic-cp} includes a rich class of problems as special cases. For example, when $\mathcal{K} = \Re_{+}^{m_1}\times\{0\}^{m_2}$ for some $m_1$ and $m_2$, $g(x) = (g_1(x),\ldots,g_{m_1}(x), h_1(x),\ldots,h_{m_2}(x))^T$ with convex $g_i$'s and affine $h_j$'s, and $P(x)$ is the indicator function of a simple convex set $X\subseteq\Re^n$, problem \eqref{eq:conic-cp} reduces to an ordinary convex programming problem
\[
\begin{array}{rl}
\min\limits_{x \in X} \{f(x): g_i(x) \leq 0, \ i=1,\ldots,m_1; h_j(x) = 0, \ j=1,\ldots,m_2\}.
\end{array}
\]

Augmented Lagrangian (AL) methods have been widely regarded as a highly effective method for solving constrained nonlinear programming (e.g., see \cite{Bertsekas99,ruszczynski2006nonlinear,nocedal-wright}). The classical AL method was initially proposed by Hestenes \cite{hestenes1969multiplier} and Powell \cite{powell1969method}, and has been extensively studied in the literature (e.g., see \cite{Rock76b,Bertsekas82}). Recently, AL methods have been applied to solve some instances of  problem \eqref{eq:conic-cp} arising in various applications such as  image processing \cite{chan2011augmented} and optimal control \cite{ito1990augmented}. 
They have also been used to solve large-scale conic programming problems (e.g., see \cite{burer2003nonlinear,jarre2008augmented,zhao2010newton}).

When applied to problem \eqref{eq:conic-cp}, AL methods proceed in the following manner. Let $\{\rho_k\}$ be a sequence of nondecreasing positive scalars and $\lambda^0\in\mcK^*$ an initial guess of the Lagrangian multiplier of \eqref{eq:conic-cp}. At the $k$th iteration, $x^{k+1}$ is obtained by approximately solving the AL subproblem
\beq \label{AL-subprob} 
\min_x\mcL(x,\lambda^k;\rho_k),
\eeq 
where $\mcL(x,\lambda;\rho)$ is the AL function of \eqref{eq:conic-cp} defined as (e.g., see \cite[Section 11.K]{rockafellar2009variational} and \cite{Shapiro04})
\[
\mcL(x,\lambda;\rho) := f(x) + P(x) + \frac{1}{2\rho}\Big[\dist^2\big(\lambda + \rho g(x), -\mcK\big) - \|\lambda\|^2\Big],
\]
and $\dist(z, -\mcK) = \min\{\|z + x\|: x\in\mcK\}$ for any $z\in\Re^m$
Then $\lambda^{k+1}$ is updated by 
\[
\lambda^{k+1} = \Pi_{\mcK^*}(\lambda^k + \rho_kg(x^{k+1})),
\]
where $\Pi_{\mcK^*}(\cdot)$ is the projection operator onto  $\mcK^*$. The iterations that update $\{\lambda^k\}$ are commonly called the \emph{outer iterations} of AL methods, and the iterations of an iterative scheme for solving AL subproblem \eqref{AL-subprob} are referred to as the \emph{inner iterations} of AL methods. In the context of large-scale optimization, a first-order method is often used to approximately solve the AL subproblem \eqref{AL-subprob} and the resulting entire method is usually called a \emph{first-order inexact AL (I-AL) method}.

In this paper we focus on developing first-order I-AL methods and studying their iteration complexity, which is an upper bound on the total number of first-order inner iterations for finding an $\epsilon$-Karush-Kuhn-Tucker ($\epsilon$-KKT) solution of  \eqref{eq:conic-cp}, that is, a primal-dual solution $(x,\lambda)$ satisfying
\begin{equation}
\label{eq:eps-approx-kkt-point}
\dist(0, \nabla f(x) + \partial P(x) + \nabla g(x)\lambda) \leq \epsilon, \quad \dist(g(x), \mcN_{\mcK^*}(\lambda))\leq \epsilon, \quad (x,\lambda) \in \dom(P) \times \mcK^*
\end{equation}
for some prescribed tolerance $\epsilon>0$. The condition \eqref{eq:eps-approx-kkt-point} is often checkable in practice and has been broadly used as a termination criterion for I-AL type of methods (e.g., see~\cite{zhao2010newton}).  It has been shown that under some mild error bound condition any point $x$ satisfying \eqref{eq:eps-approx-kkt-point} with a small $\epsilon$ is close to an optimal solution of problem \eqref{eq:conic-cp} (e.g., see~\cite{Rock76b}). As problem \eqref{eq:conic-cp} arising in various applications is of large scale, a first-order I-AL method with a low complexity bound for finding an $\epsilon$-KKT solution of  \eqref{eq:conic-cp} is highly desirable. The main contributions of this paper consist of: (i) proposing a first-order I-AL method that resembles the classical AL method and establishing its iteration complexity, which reveals how good iteration complexity of the classical I-AL method can be; (ii) proposing an adaptively regularized first-order I-AL method that achieves a significantly improved iteration complexity over existing first-order I-AL methods for finding an $\epsilon$-KKT solution of problem \eqref{eq:conic-cp}; and (iii) a technically new complexity analysis of the I-AL methods based on a sharp analysis of inexact proximal point algorithm (PPA) and the connection between the I-AL methods and inexact PPA, which provides more insights than existing complexity analyses of the first-order I-AL methods in the literature that typically regard the I-AL methods as an inexact dual gradient method.

\subsection{Related works}

Aybat and Iyengar \cite{aybat2013augmented} proposed a first-order I-AL method for solving a special case of \eqref{eq:conic-cp} with affine mapping $g$. 
In particular, they applied an optimal first-order method (e.g., see \cite{Nesterov05,Tseng08}) to find an 
approximate solution $x^{k+1}$ of the AL subproblem \eqref{AL-subprob} such that 
\[
\mcL(x^{k+1},\lambda^k;\rho_k) - \min_x\mcL(x,\lambda^k;\rho_k) \leq \eta_k
\] 
for some $\eta_k>0$. It is shown in \cite{aybat2013augmented} that this method with some suitable 
choice of $\{\rho_k\}$ and $\{\eta_k\}$ can find an approximate solution $x$ of \eqref{eq:conic-cp} satisfying 
\begin{equation}
\label{eq:eps-subopt}
|F(x) - F^*| \leq \epsilon, \quad \dist(g(x),-\mcK) \leq \epsilon 
\end{equation}
for some $\epsilon>0$ in at most $\mathcal{O}(\epsilon^{-1}\log\epsilon^{-1})$ first-order inner iterations. In addition, Necoara et al.\ \cite{necoara2015complexity} proposed an accelerated first-order I-AL method for solving the same problem as considered in \cite{aybat2013augmented}, in which an acceleration  scheme \cite{DeGlNe14} is applied to $\{\lambda^k\}$ for possibly better convergence. It is claimed in 
 \cite{necoara2015complexity} that this method with a suitable choice of $\{\rho_k\}$ and $\{\eta_k\}$ can find an approximate solution $x$ of \eqref{eq:conic-cp} satisfying \eqref{eq:eps-subopt} 
in at most $\mathcal{O}(\epsilon^{-1})$ first-order inner iterations. More recently, Xu \cite{Xu17} proposed an I-AL method for solving a special case of \eqref{eq:conic-cp} with $\mcK$ being the nonnegative orthant, which can find an approximate solution $x$ of \eqref{eq:conic-cp} satisfying \eqref{eq:eps-subopt} in at most $\mathcal{O}(\epsilon^{-1})$ first-order inner iterations. Some other related works on I-AL type of methods can be found, for example, in \cite{LLM17,PNT17,Xu17b}.

Since $F^*$ is typically unknown, the criterion \eqref{eq:eps-subopt} is not checkable and cannot be used as a practical termination criterion for the I-AL methods \cite{aybat2013augmented,necoara2015complexity,PNT17,Xu17b,Xu17} in general. Thus, for the practical use of these methods, one has to terminate them by a checkable criterion, which may result in a substantially different solution from an $\epsilon$-optimal solution defined in \eqref{eq:eps-subopt}. Due to this, it may be challenging for them to find such an $\epsilon$-optimal solution \emph{in practice}. Moreover, for these I-AL methods, $\{\rho_k\}$ and $\{\eta_k\}$ are specifically chosen to achieve a low first-order iteration complexity with respect to \eqref{eq:eps-subopt}. Such a choice may however not lead to a low first-order iteration complexity with respect to a checkable termination criterion. Therefore, the iteration-complexity results obtained in \cite{aybat2013augmented,necoara2015complexity,PNT17,Xu17b,Xu17} 
with respect to the criterion \eqref{eq:eps-subopt} do not seem to have much practical merits in general.

In addition to the aforementioned I-AL methods, Lan and Monteiro \cite{Lan16} proposed a first-order I-AL method for finding an \emph{$\epsilon$-KKT solution} of a special case of \eqref{eq:conic-cp} with $g=\cA(\cdot)$, $\mcK = \{0\}^m$ and $P$ being the indicator function of a simple compact convex set $X$, that is, 
\begin{equation}\label{eq:conic-lp}
\min\left\{f(x): \cA(x)= 0, x\in X\right\},
\end{equation}
where $\cA:\Re^n \to \Re^m$ is an affine mapping. Roughly speaking, their I-AL method consists of two stages, particularly, primary stage and postprocessing stage.  The primary stage is to execute the usual I-AL steps similar to those in \cite{aybat2013augmented} 
but with static $\rho_k \equiv \cO(D_\Lambda^{3/4}\epsilon^{-3/4})$ and $\eta_k \equiv  \cO(D_\Lambda^{1/4}\epsilon^{7/4})$ 
until a certain approximate solution $(\tx,\tlambda)$ is found,\footnote{It means that $\rho_k=\rho$ and $\eta_k=\eta$  for all $k$ for some $\rho=\cO(D_\Lambda^{3/4}\epsilon^{-3/4})$ and $\eta = \cO(D_\Lambda^{1/4}\epsilon^{7/4})$.} where $D_\Lambda = \min\{\|\lambda^0 - \lambda\|: \lambda\in\Lambda^*\}$ and $\Lambda^*$ is the set of optimal solutions of the Lagrangian dual problem associated with problem \eqref{eq:conic-lp}. The postprocessing stage is mainly to execute a single I-AL step with 
a penalty parameter $\rho=\cO(D_\Lambda^{3/4}\epsilon^{-3/4})$ and an AL subproblem tolerance parameter $\eta= \cO(\min(D_\Lambda^{3/4}\epsilon^{5/4}, D_\Lambda^{-3/4}\epsilon^{11/4}))$, starting with $(\tx,\tlambda)$.  
It is shown in \cite{Lan16} that this I-AL method can find  an \emph{$\epsilon$-KKT solution} of \eqref{eq:conic-lp} in at most $\mathcal{O}(\epsilon^{-7/4})$ first-order inner iterations in theory. However, this method is much less practical than the classical I-AL method. Indeed, in the primary stage, this I-AL method uses $\rho_k$ and $\eta_k$ of same value through all outer iterations, which may be respectively overly large and small. Such a choice of $\rho_k$ and $\eta_k$ is clearly against the common practical choice that $\rho_0$ and $\eta_0$ are relatively small and large, respectively, and $\rho_k$ gradually increases and $\eta_k$ gradually decreases as iteration progresses. Moreover,  $\rho_k$ and $\eta_k$ in this method require some knowledge of $D_\Lambda$, which is not known a priori and needs to be estimated by a sophisticated and expensive ``guess-and-check'' procedure proposed in \cite{Lan16}. These two aspects evidently make this I-AL method much more sophisticated and less practical than the classical I-AL method.

Besides, Lan and Monteiro \cite{Lan16} proposed a modified I-AL method by applying their 
aforementioned first-order I-AL method with $D_\Lambda$ replaced by $D^\epsilon_\Lambda$ to the perturbed problem 
\begin{equation}\label{eq:conic-lp-pert}
\min\left\{f(x)+\frac{\epsilon}{4D_X}\|x- x^0\|^2: \cA(x)= 0, x\in X\right\},
\end{equation}
starting with some $(x^0,\lambda^0)$,
where $D_X=\max\{\|x-y\|: x,y\in X\}$, $D^\epsilon_\Lambda=\min\{\|\lambda^0-\lambda\|:\lambda\in\Lambda^*_\epsilon\}$, and $\Lambda^*_\epsilon$ is the set of Lagrangian dual optimal solutions associated with problem \eqref{eq:conic-lp-pert}.  They showed that their modified I-AL method can find  an \emph{$\epsilon$-KKT solution} of \eqref{eq:conic-lp} in at most 
\beq \label{complexity-mAL}
\mathcal{O}\left\{\left(\frac{\sqrt{D^\epsilon_\Lambda}}{\epsilon} \left[\log\frac{\sqrt{D^\epsilon_\Lambda}}{\epsilon}\right]^{\frac34}+\frac{1}{\sqrt{\epsilon}}\log\frac{\sqrt{D^\epsilon_\Lambda}}{\epsilon}\right) \max\left(1,\log\log\frac{\sqrt{D^\epsilon_\Lambda}}{\epsilon}\right)\right\}
\eeq
first-order inner iterations. Since $D^\epsilon_\Lambda$ depends on $\epsilon$ (see an example in Appendix \ref{app-example}) and its order dependence on $\epsilon$ is generally unknown, it is not clear about 
the order dependence of the iteration complexity \eqref{complexity-mAL} on $\epsilon$. 

\subsection{Main contribution}

The goal of this paper is to propose first-order I-AL methods with significantly improved iteration complexity over existing first-order I-AL methods for finding an \emph{$\epsilon$-KKT solution} of problem \eqref{eq:conic-cp}.  Our main contribution is listed below.

\bi
\item We propose a first-order I-AL method that resembles the classical I-AL method establish its first-order iteration complexity $\mathcal{O}(\epsilon^{-7/4})$ for finding an \emph{$\epsilon$-KKT solution} of problem \eqref{eq:conic-cp}. Due to the similarity between this I-AL method and the classical I-AL method,  our complexity result reveals that the iteration complexity of the classical I-AL method for finding an $\epsilon$-KKT solution of \eqref{eq:conic-cp} appears to be $\mathcal{O}(\epsilon^{-7/4})$.

\item We propose an \emph{adaptively regularized} first-order I-AL method that establish its first-order iteration complexity $\mathcal{O}(\epsilon^{-1}\log \epsilon^{-1})$ for finding an \emph{$\epsilon$-KKT solution} of problem \eqref{eq:conic-cp}, which significantly improves the previously best-known iteration-complexity $\mathcal{O}(\epsilon^{-7/4})$ achieved by  first-order I-AL methods for finding an $\epsilon$-KKT solution of  \eqref{eq:conic-cp}. This complexity result implies that the adaptively regularized first-order I-AL method is generally superior to the classical I-AL method for finding an $\epsilon$-KKT solution of \eqref{eq:conic-cp}.

\item Our complexity analysis of the I-AL methods is technically new, which is based on a sharp analysis of inexact PPA and the connection between the I-AL methods and inexact PPA. It is vastly different from existing complexity analyses of the first-order I-AL methods, which typically regard the I-AL methods as an inexact dual gradient method. Since the operator associated with the monotone inclusion problem linked to the I-AL methods is closely related to the KKT conditions, our analysis is more appropriate and provides more insights than existing ones in the literature.
\ei

\subsection{Outline}
The rest of this paper is organized as follows. In Section \ref{sec:preliminaries} we introduce some notation and the concept of  $\epsilon$-KKT solution. In Section \ref{sec:i-al} we propose a first-order I-AL method and present its iteration complexity. Also, in Section \ref{sec:mi-al} we propose an adaptively regularized first-order I-AL method and present its iteration complexity. In Section \ref{sec:proof} we provide a proof for the technical results stated in Sections \ref{sec:i-al} and \ref{sec:mi-al}. In Section \ref{sec:numerical} we present some numerical results for the proposed algorithms. Finally, we make some concluding remarks in Section \ref{sec:conclude}.

\section{Notation and preliminaries}\label{sec:preliminaries}
The following notations will be used throughout this paper. Let $\Re^n$ denote the Euclidean space of dimension $n$, $\langle\cdot,\cdot\rangle$ denote the standard inner product, and $\|\cdot\|$ stand for the Euclidean norm. The symbols $\Re_+$ and $\Re_{++}$ stand for the set of nonnegative and positive real numbers, respectively.  

Given a closed convex function $h:\Re^n\rightarrow(-\infty,\infty]$, $\partial h$ and $\dom (h)$ denote the subdifferential and domain  of $h$, respectively.  The proximal operator associated with $h$ is denoted by  
$\prox_h$,  that is,
\begin{equation*}
\label{eq:def-prox}
\prox_h(z) = \argmin_{x\in\Re^n} \left\{ \frac{1}{2}\|x - z\|^2 + h(x) \right\}, \quad \forall z \in \Re^n.
\end{equation*}
Given a non-empty closed convex set $C\subseteq\Re^n$, $\dist(z,C)$ stands for the Euclidean distance from $z$ to $C$, and $\Pi_C(z)$ denotes the Euclidean projection of $z$ onto $C$, namely, 
\[
 \Pi_C(z) = \argmin\{\|z - x\|: x\in C\}, \quad \dist(z,C)=\left\|z - \Pi_C(z)\right\|, \quad \forall z \in \Re^n.
\]
The normal cone of $C$ at any $z\in C$ is denoted by $\mcN_C(z)$. For the closed convex cone $\mcK$, we use $\mcK^*$ to denote the dual cone of $\mcK$, that is,  $\mcK^* = \{y\in\Re^m: \langle y,x\rangle\geq 0, \; \forall x\in\mcK\}$.

The Lagrangian function $l(x,\lambda)$ of problem~\eqref{eq:conic-cp} is defined as
\begin{equation*}
l(x,\lambda) = \left\{
\begin{array}{ll}
f(x) + P(x) + \Inner{\lambda}{g(x)} & \mbox{if}\; x\in \dom(P) \mbox{ and } \lambda \in \mcK^{*}, \\
-\infty & \mbox{if}\; x\in \dom(P) \mbox{ and } \lambda\notin\mcK^{*}, \\
\infty & \mbox{if}\; x\notin \dom(P),
\end{array}
\right.
\end{equation*}
which is a closed convex-concave function. The Lagrangian dual function $d:\Re^m\rightarrow[-\infty,\infty)$ is defined as 
\begin{equation}
\label{eq:lag-dual-func}
d(\lambda) = \inf_x l(x,\lambda) = \left\{
\begin{array}{ll}
\inf\limits_x\  \{ f(x) + P(x) + \Inner{\lambda}{g(x)}\} & \mbox{if}\; \lambda\in\mcK^{*},\\
-\infty & \mbox{if}\; \lambda\notin\mcK^{*},
\end{array}
\right.
\end{equation}
which is a closed concave function. Moreover, the augmented Lagrangian function for problem~\eqref{eq:conic-cp} is defined as (e.g., see~\cite{Shapiro04})
\begin{equation}
\label{eq:al-func}
\mcL(x,\lambda;\rho) = f(x) + P(x) + \frac{1}{2\rho}\Big[\dist^2\big(\lambda + \rho g(x), -\mcK\big) - \|\lambda\|^2\Big],
\end{equation}
where $\rho>0$ is  a penalty parameter. For convenience, we let
\beq \label{S-fun}
\mcS(x,\lambda;\rho):= f(x) + \frac{1}{2\rho} \dist^2\big(\lambda+\rho g(x), -\mcK \big).
\eeq
It is clear to see that 
\[
\mcL(x,\lambda;\rho) = \mcS(x,\lambda;\rho) + P(x) - \frac{\|\lambda\|^2}{2\rho}.
\]
The augmented Lagrangian dual function of~\eqref{eq:conic-cp} is given by
\begin{equation}
\label{eq:aug-dual-func}
d(\lambda;\rho)= \min_{x\in\Re^n} \mcL(x,\lambda;\rho).
\end{equation}
The Lagrangian dual problem \eqref{dual-prob} can thus be rewritten as  
\begin{equation}\label{eq:dual-conic-cp}
d^{*} = \max_\lambda\; d(\lambda).
\end{equation}
Let $\partial l:\Re^n\times\Re^m\rightrightarrows\Re^n\times \Re^m$ and $\partial d:\Re^m\rightrightarrows\Re^m$ be respectively the subdifferential mappings associated with $l$ and $d$ (e.g., see \cite{R97}). We define two set-valued operators associated with problems \eqref{eq:conic-cp} and \eqref{dual-prob} as follows:
\beqa
&&\mcT_d: \lambda \rightarrow \{ u\in\Re^m: -u\in \partial d(\lambda)\},  \quad \forall \lambda\in\Re^m, \label{eq:def-two-opt1} \\
&& \mcT_l: (x,\lambda) \rightarrow \{ (v,u)\in\Re^n\times\Re^m: (v,-u)\in\partial l(x,\lambda)\}, \quad \forall (x,\lambda) \in \Re^n \times \Re^m. \label{eq:def-two-opt2}
\eeqa
It is well known that $\lambda^*$ is an optimal solution of the Lagrangian dual problem~\eqref{eq:dual-conic-cp} if and only if $0\in \partial d(\lambda^*)$, and $(x^*,\lambda^*)$ is a saddle point\footnote{$(x^*,\lambda^*)$ is called a saddle point of $l$ if it satisfies $ \sup_{\lambda}l(x^*,\lambda) = l(x^*,\lambda^*) = \inf_x l(x,\lambda^*)$.} of $l$ if and only if $(0,0)\in\partial l(x^*,\lambda^*)$. In addition, it can be verified that 
\begin{equation}
\label{eq:express-Tl}
\partial l(x,\lambda) = \left\{
\begin{array}{ll}
\left(
\begin{array}{c}
\nabla f(x) + \partial P(x) + \nabla g(x)\lambda \\
g(x) - \mcN_{\mcK^*}(\lambda)
\end{array}
\right), & \mbox{if}\; x\in\dom(P)\; \mbox{and}\; \lambda\in\mcK^*, \medskip
\\ 
\emptyset, & \mbox{otherwise}.
\end{array}\right.
\end{equation}
This enables us to write the KKT optimality condition for problem \eqref{eq:conic-cp} as follows.
\begin{prop}
Under Assumption~\ref{ass:function}, $x^*\in\Re^n$ is an optimal solution of~\eqref{eq:conic-cp} if and only if there exists $\lambda^*\in\Re^m$ such that
\begin{equation}
\label{eq:KKT-condition}
(0,0)\in \partial l(x^*,\lambda^*),
\end{equation}
or equivalently, $(x^*,\lambda^*)$ satisfies the KKT conditions for~\eqref{eq:conic-cp}, that is,
\[
0\in\nabla f(x^*) + \partial P(x^*) + \nabla g(x^*) \lambda^*, \quad \lambda^* \in\mcK^*, \quad g(x^*)\preceq_\mcK 0, \quad \langle \lambda^*, g(x^*) \rangle = 0.
\]
\end{prop}
\begin{proof}
The result \eqref{eq:KKT-condition} follows from \cite[Theorem 36.6]{R97}. By \eqref{eq:express-Tl}, it is not hard to see that \eqref{eq:KKT-condition} holds if and only if $0\in \nabla f(x^*) + \partial P(x^*) + \nabla g(x^*)\lambda^*$, $\lambda^*\in\mcK^*$, and $g(x^*) \in\mcN_{\mcK^*}(\lambda^*)$. By the definition of $\mcK^*$ and $\mcN_{\mcK^*}$, one can verify that $g(x^*)\in\mcN_{\mcK^*}(\lambda^*)$ is equivalent to $g(x^*)\preceq_\mcK 0$ and $\langle\lambda^*, g(x^*) \rangle = 0$. The proof is then completed.
\end{proof}
 
In practice it is generally impossible to find an exact KKT solution $(x^*,\lambda^*)$ satisfying \eqref{eq:KKT-condition}.  We are instead interested in seeking an approximate KKT solution 
of~\eqref{eq:conic-cp} that is defined as follows.

\begin{defi}
	\label{defi:approx-KKT}
	Given any $\epsilon>0$, we say $(x,\lambda)\in\Re^n\times\Re^m$ is an $\epsilon$-KKT solution of~\eqref{eq:conic-cp}  if there exists $(u,v)\in\partial l(x,\lambda)$ such that $\|u\|\leq \epsilon$ and $\|v\|\leq \epsilon$.
\end{defi}

\begin{rem} \label{rem1} 
	\begin{itemize}
		\item[(a)] By \eqref{eq:express-Tl} and Definition \ref{defi:approx-KKT}, one can see that $(x,\lambda)$ is an $\epsilon$-KKT solution of~\eqref{eq:conic-cp}  if and only if $x\in\dom(P)$, $\lambda\in\mcK^*$, $\dist(0, \nabla f(x) + \partial P(x) + \nabla g(x)\lambda) \leq \epsilon$, and $\dist(g(x), \mcN_{\mcK^*}(\lambda))\leq \epsilon$.  It reduces 
		to an $\epsilon$-KKT solution introduced in~\cite{Lan16} when $g$ is affine and $\mcK = \{0\}$.
		\item[(b)] For a given $(x,\lambda)$, it is generally not hard to verify whether it is an $\epsilon$-KKT solution of~\eqref{eq:conic-cp}. Therefore, Definition \ref{defi:approx-KKT} gives rise to a checkable termination criterion \eqref{eq:eps-approx-kkt-point} that will be used in this paper.
	\end{itemize}
\end{rem}

\section{A first-order I-AL method and its iteration complexity}\label{sec:i-al}

In this section we propose a first-order I-AL method that resembles the classical I-AL method, and study its first-order iteration complexity for finding an $\epsilon$-KKT solution of problem~\eqref{eq:conic-cp}.

Recall Remark \ref{rem1}(a) that $(x,\lambda)$ is an $\epsilon$-KKT solution of~\eqref{eq:conic-cp} if and only if it satisfies that $x\in\dom(P)$, $\lambda\in\mcK^*$, $\dist(g(x), \mcN_{\mcK^*}(\lambda))\leq \epsilon$, and $\dist(0, \nabla f(x) + \partial P(x) + \nabla g(x)\lambda) \leq \epsilon$. In what follows, we propose an I-AL method to 
find a pair $(x,\lambda)$ satisfying these conditions. Given that the proximal operator 
associated with $P$ and the projection onto $\mcK$ can be exactly evaluated (see Assumption \ref{ass:function}), the first two conditions can be easily satisfied by the iterates of our I-AL method. Observe that the last condition is generally harder to satisfy than the third one since it involves $\nabla f$, $\nabla g$ and $\partial P$. Due to this, our I-AL method consists of two stages, particularly, primary stage and postprocessing stage. In the primary stage, the AL subproblems are solved roughly and the usual I-AL steps are executed until either an $\epsilon$-KKT solution of \eqref{eq:conic-cp} is obtained or a pair $(x^k,\lambda^k)$ satisfying nearly the third condition but roughly the last condition is found. In the postprocessing stage,  the last AL subproblem arising in the primary stage is resolved to a higher accuracy for obtaining some point $\tilde x$, starting with $x^k$, and a proximal step is 
then applied to $\mcL(\cdot,\lambda^k,\rho_k)$ at $\tilde x$ and to $l(\tilde x, \cdot)$ at $\lambda^k$ respectively, to generate the output $(x^+,\lambda^+)$. 
 
Our first-order I-AL method for solving problem \eqref{eq:conic-cp} is presented as follows.

\begin{algo}[A first-order I-AL method]
\label{alg:i-al}
\normalfont
\mbox{}
\begin{itemize}
	\item[0.] Input $\epsilon>0$, $\lambda^0\in\mcK^*$, nondecreasing $\{\rho_k\} \subset \Re_{++}$, and $0<\eta_k \downarrow 0$. Set $k=0$.

	\item[1.] Apply Algorithm~\ref{alg:opt-fom} (in Appendix \ref{app-fom}) to the problem $\min_{x} \mcL(x,\lambda^k;\rho_k)$ to find $x^{k+1}\in\dom(P)$ satisfying
	\begin{equation}
	\label{eq:inext-subprob}
	\mcL(x^{k+1},\lambda^k;\rho_k) - \min_{x} \mcL(x,\lambda^k;\rho_k) \leq \eta_k.\footnote{In view of Proposition \ref{prop:cplx-fom-coro},  one can terminate Algorithm~\ref{alg:opt-fom} once a point $x^{k+1}$ satisfying $\Psi(x^{k+1})-\underline\Psi_{k+1} \le \eta_k$ is found, where $\Psi(\cdot)=\mcL(\cdot,\lambda^k;\rho_k)$, and $\underline\Psi_{k+1}$ is defined according to \eqref{eq:cplx-bd-fom} that can be cheaply computed (see Remark \ref{phik}). Such $x^{k+1}$ clearly satisfies \eqref{eq:inext-subprob}.}
	\end{equation}
	\item[2.] Set $\lambda^{k+1} = \Pi_{\mcK^*}(\lambda^k + \rho_kg(x^{k+1}))$. 
       \item[3.] If $(x^{k+1},\lambda^{k+1})$ satisfies \eqref{eq:eps-approx-kkt-point}, output $(x^+,\lambda^+)=(x^{k+1},\lambda^{k+1})$ and terminate.                    
	\item[4.] If the following inequalities are satisfied
	\begin{equation}
	\label{eq:check-stop}
	\frac{1}{\rho_k}\|\lambda^{k+1} - \lambda^k\| \leq \frac{3}{4}\epsilon, \quad\quad  \frac{\eta_k}{\rho_k}\leq \frac{\epsilon^2}{128},
	\end{equation}
call the subroutine $(x^+,\lambda^+) =$ Postprocessing$(\lambda^k,\rho_k,x^{k+1},\epsilon)$, output $(x^+,\lambda^+)$ and terminate.
	\item[5.] Set $k\leftarrow k+1$ and go to Step 1.
\end{itemize}
{\bf End.}
\end{algo}

We next present the subroutine Postprocessing that is used in Step 4 of Algorithm \ref{alg:i-al}. Before proceeding, let $L_g$ be the  Lipschitz constant of $g$ on $\dom(P)$ and $M_g = \max\limits_{x\in\dom(P)}\|g(x)\|$.  

\medskip\noindent
Subroutine $(x^+,\lambda^+) =$ {\bf Postprocessing}$(\tilde{\lambda},\tilde{\rho},\tilde{x},\epsilon)$
\begin{itemize}
        \item[0.] Input $\tilde{\lambda}\in \mcK^*$, $\tilde{\rho}>0$, $\tilde{x}\in\dom(P)$, and $\epsilon>0$.
       \item[1.] Set 
\begin{equation}
	\label{eq:postpro-para}
\tilde{L} = L_{\nabla f} + L_{\nabla g}\big(\|\tilde{\lambda}\| + \tilde{\rho}M_g\big) + \tilde{\rho}L_g^2, \quad \quad\tilde{\eta} = \epsilon^2 \cdot \min \left\{ \frac{\tilde{\rho}}{128}, \frac{1}{8\tilde{L}}\right\}.
	\end{equation}
	\item[2.] Apply Algorithm~\ref{alg:opt-fom} to the problem $\min_{x} \mcL(x, \tilde{\lambda}; \tilde{\rho})$ starting with $\tilde{x}$ to find $\hat{x}$ such that 
	\begin{equation}
	\label{eq:inext-postpro}
	\mcL(\hat{x}, \tilde{\lambda}; \tilde{\rho})- \min_{x} \mcL(x, \tilde{\lambda}; \tilde{\rho})\leq \tilde{\eta}.
	\end{equation}
	\item[3.] Output the pair $(x^+,\lambda^+)$, which is computed by
	\begin{equation}
	\label{eq:output-postpro}
	x^+ = \prox_{P/\tilde{L}}\big(\hat{x} - \nabla_x\mcS(\hat{x},\tilde{\lambda};\tilde{\rho})/\tilde{L}\big), \quad \lambda^+ = \Pi_{\mcK^*}\big(\tilde{\lambda} + \tilde{\rho} g(x^+) \big),
	\end{equation}
where $\mcS$ is defined in \eqref{S-fun}.
\end{itemize}
{\bf End.}

\medskip
For ease of later reference, we refer to the first-order iterations of Algorithm~\ref{alg:opt-fom} for solving the 
AL subproblems as the {\it inner iterations} of Algorithm \ref{alg:i-al}, and call, the update from $(x^k,\lambda^k)$ to $(x^{k+1},\lambda^{k+1})$ or the postprocessing step, an {\it outer iteration} of Algorithm \ref{alg:i-al}.  We now make some remarks on Algorithm~\ref{alg:i-al} 
as follows.

\begin{rem} 
\begin{itemize}
\item[(a)] The subroutine Postprocessing is inspired by~\cite{Lan16}, in which a similar procedure is proposed for solving a special case of problem \eqref{eq:conic-cp} with affine $g$
 and $\mcK = \{0\}$. The main purpose of this subroutine is to obtain a better iteration  complexity.
	\item[(b)] Compared to the I-AL method in \cite{Lan16}, our I-AL method is much simpler and more closely resembles the classical I-AL method. Indeed, the I-AL method \cite{Lan16} uses the static  $\rho_k \equiv \cO(D_\Lambda^{3/4}\epsilon^{-3/4})$ and $\eta_k \equiv \cO(D_\Lambda^{1/4}\epsilon^{7/4})$  through all outer iterations in the primary stage,\footnote{It means that $\rho_k=\rho$ and $\eta_k=\eta$  for all $k$ for some $\rho=\cO(D_\Lambda^{3/4}\epsilon^{-3/4})$ and $\eta = \cO(D_\Lambda^{1/4}\epsilon^{7/4})$.} where $D_\Lambda = \min\{\|\lambda^0 - \lambda\|: \lambda\in\Lambda^*\}$ and $\Lambda^*$ is the set of optimal solutions of the Lagrangian dual problem associated with problem \eqref{eq:conic-lp}. Such $\{\rho_k\}$ and $\{\eta_k\}$ may be overly large and small, respectively. This is clearly against the common practical choice that $\rho_0$ and $\eta_0$ are relatively small and large, respectively, and $\{\rho_k\}$ gradually increases and $\{\eta_k\}$ progressively decreases. Moreover, the above choice of $\rho_k$ and $\eta_k$ requires some knowledge of $D_\Lambda$, which is not known a priori and needs to be estimated by
a sophisticated and expensive ``guess-and-check'' procedure proposed in \cite{Lan16}.
In contrast, our I-AL method uses a practical choice of $\{\rho_k\}$ and $\{\eta_k\}$, which dynamically change throughout the iterations. Also, it does not use any knowledge of $D_\Lambda$ and thus a ``guess-and-check'' procedure is not required.
\end{itemize}
\end{rem}

In the rest of this section we present our main results for Algorithm \ref{alg:i-al}, whose proof is deferred  to Subsection \ref{sec:pf-i-al}. To proceed, we introduce some further notation below. Let $\Lambda^*$ be the set of optimal solutions of problem \eqref{dual-prob} and $\hat\lambda^*\in \Lambda^*$ such that $\|\lambda^0-\hat\lambda^*\|=\dist(\lambda^0,\Lambda^*)$. In addition, we define 
\begin{align}
 & D_X = \max\limits_{x,y\in\dom(P)}\|x - y\|,  \  D_\Lambda = \|\lambda^0 - \hat\lambda^*\|, \ B = L_{\nabla f} + L_{\nabla g}(\|\hat\lambda^*\| + D_\Lambda), \label{param} \\
 & C = L_{\nabla g}M_g + L_g^2, \ \bD=\max\{D_\Lambda,1\}, \ \bB =\max\{B,1\}, \ \bC =\max\{C,1\}, \label{bparam}
\end{align}
where  $L_{\nabla f}$, $L_{\nabla g}$, $L_g$ and $M_g$ are defined above.
We start with the following theorem which shows that Algorithm \ref{alg:i-al} with a suitable choice of $\{\rho_k\}$ and $\{\eta_k\}$ is guaranteed to find an $\epsilon$-KKT solution of problem \eqref{eq:conic-cp} in a finite number of outer iterations.

\begin{thm}
	\label{prop:output-gantee}
	\bi
	\item[(i)]
	If  Algorithm~\ref{alg:i-al} successfully terminates, then the output $(x^+,\lambda^+)$  is an $\epsilon$-KKT solution of problem \eqref{eq:conic-cp}.
	\item[(ii)]
	Suppose that $\{\rho_k\}$  and $\{\eta_k\}$ satisfy that
	\beq \label{rho-eta}
	\rho_k>0 \ \mbox{is nondecreasing}, \quad 0< \frac{\eta_k}{\rho_k} \to 0, \quad 
	\frac{\sum_{i=0}^{2k}\sqrt{\rho_i\eta_i}}{\rho_k\sqrt{k+1}} \to 0.  \footnote{For example,  $\rho_k=\hat{C} (k+1)^{3/2}$ and $\eta_k={\tilde C}(k+1)^{-5/2}$ satisfy \eqref{rho-eta} for any $\hat{C} $, ${\tilde C}>0$.}
	\eeq 
	Then Algorithm~\ref{alg:i-al} terminates  in a finite number of outer iterations. 
	\item[(iii)] Furthermore, if $N$  is a nonnegative integer such that
	\begin{equation}
	\label{eq:cplx-out-general}
	\frac{D_\Lambda + 2\sum_{k=0}^{2N}\sqrt{2\rho_k\eta_k}}{\rho_N\sqrt{N+1}} \leq \frac{\epsilon}{2}, \quad \frac{\eta_N}{\rho_N}\leq \frac{\epsilon^2}{128},
	\end{equation}
then Algorithm~\ref{alg:i-al} terminates in at most $2N+1$ outer iterations. 
	\ei
\end{thm}

The next theorem provides an upper bound on the total number of the inner iterations Algorithm~\ref{alg:i-al}, that is, the total iterations of Algorithm \ref{alg:opt-fom} applied to solve all AL subproblems of Algorithm~\ref{alg:i-al}. 

\begin{thm}
	\label{thm:cplx-i-al}
	Let $\epsilon>0$ be given, and $\bC$, $D_X$, and $\bD$ be defined in \eqref{param} and 
\eqref{bparam}.  Suppose that $\{\rho_k\}$ and $\{\eta_k\}$ are chosen as
	\begin{equation}
	\label{eq:para-i-al}
	\rho_k = \rho_0(k+1)^{\frac{3}{2}}, \ \ \eta_k = \eta_0(k+1)^{-\frac{5}{2}}\cdot\min\{1,\sqrt{\epsilon}\}
	\end{equation}
	for some $\rho_0\geq 1$ and $0<\eta_0\leq 1$. Then, the total number of inner iterations of Algorithm \ref{alg:i-al} for finding an $\epsilon$-KKT solution of problem \eqref{eq:conic-cp} is at most $\mathcal{O}\big({\cal T}\big(\min\{1,\epsilon\}\big)\big)$, where
\[
{\cal T}(t) = \frac{D_X\bD^{\frac{3}{2}}\bC}{t^\frac{7}{4}} + \frac{D_X\bD^{\frac{5}{4}}\bB^{\frac{1}{2}}(1+L_{\nabla g}^{\frac{1}{2}})}{t^{\frac{11}{8}}} + \frac{D_X\bD^{\frac{1}{4}}(L_{\nabla g} + L_{\nabla g}^{\frac{1}{2}})}{t^{\frac{9}{8}}} + \frac{D_X\bB}{t} + \frac{\bD^{\frac{1}{2}}}{t^{\frac{1}{2}}}. 
\]
\end{thm}

\begin{rem}
	\bi
	\item[(i)]  
One can observe from Lemma \ref{lemma:lip-const-bd} and the proof of Theorem \ref{thm:cplx-i-al} that the worst-case upper bound of the total number of  inner iterations of Algorithm \ref{alg:i-al} is given by
\[
\sum^{2N+1}_{k=0} \ceil{D_X\sqrt{\frac{2(C\rho_k + B + L_{\nabla g}\sum_{i=0}^{k-1}\sqrt{2\rho_i\eta_j})}{\eta_k}}\ }
\]
It can be shown that $\rho_k = \cO((k+1)^{3/2})$ and  $\eta_k = \cO((k+1)^{-5/2}\min\{1,\sqrt{\epsilon}))$ minimize this quantity subject to the constraints given in \eqref{eq:cplx-out-general}. 
	\item[(ii)] From Theorem \ref{thm:cplx-i-al}, one can see that for any $\epsilon\in(0,1)$, the first-order iteration complexity of Algorithm \ref{alg:i-al} for finding an $\epsilon$-KKT solution of problem \eqref{eq:conic-cp} is $\cO(\epsilon^{-7/4})$, which is in the same order as the one for the I-AL method \cite{Lan16}. Nevertheless, Algorithm \ref{alg:i-al} is much more efficient than the latter method as observed in our numerical experiment. The main reason for this  is perhaps that Algorithm \ref{alg:i-al} uses the dynamic $\{\rho_k\}$ and $\{\eta_k\}$, while I-AL method \cite{Lan16} uses the static ones through all iterations and also needs a ``guess-and-check'' procedure to approximate the unknown parameter $D_\Lambda$.
	\ei
\end{rem}

Finally, the following theorem shows that for the $\epsilon$-KKT solution $(x^+,\lambda^+)$ outputted by Algorithm \ref{alg:i-al}, its primal objective gap, dual objective gap, and constraint violation are at most $\mathcal{O}(\epsilon)$.

\begin{thm}
	\label{thm:eps-funv}
	Consider the same setting as in Theorem \ref{thm:cplx-i-al}. Then, the output $(x^+,\lambda^+)$ of Algorithm \ref{alg:i-al} satisfies 
	\begin{align}
		\dist(g(x^+),-\mcK) & \leq \epsilon, \qquad -C_1\epsilon \leq F(x^+) - F^* \leq (D_X+ C_2\cdot\max\{1,\epsilon\})\epsilon, \label{eq:primal-val-bd} \\
	0  \leq F^* - d(\lambda^+) & \leq (D_X+ C_1 + C_2\cdot\max\{1,\epsilon\})\epsilon, \label{eq:dual-val-bd}
	\end{align}
	where $C_1= \dist(0, \Lambda^*)$ and $C_2 = 84\rho_0\bar{D}_\Lambda + \|\hat{\lambda}^*\|$.
\end{thm}

\section{An adaptively regularized I-AL method with improved iteration complexity}
\label{sec:mi-al}

In this section, we propose an adaptively regularized first-order I-AL method and show that it achieves a significantly improved first-order iteration complexity than Algorithm \ref{alg:i-al} and existing first-order I-AL methods in the literature for finding an $\epsilon$-KKT solution of \eqref{eq:conic-cp}. In particular, at each $k$th outer iteration it modifies Algorithm \ref{alg:i-al} by adding a regularization term $\|x - x^k\|^2/(2\rho_k)$ to the AL function $\mcL(x,\lambda^k;\rho_k)$ and also solving the AL subproblem to a higher accuracy. Moreover, it does not need a postprocessing stage. Since the regularization terms change dynamically, it is substantially different from those in \cite{necoara2015complexity,Lan16,Xu17}.

Our adaptively regularized first-order I-AL method for problem \eqref{eq:conic-cp} is presented as follows.

\begin{algo}[An adaptively regularized I-AL method]
	\label{alg:mi-al}
	\normalfont
	\mbox{}
	\begin{itemize}
		\item[0.] Input $\epsilon>0$, $(x^0,\lambda^0)\in\dom(P)\times \mcK^*$,  nondecreasing $\{\rho_k\} \subset \Re_{++}$, and $0<\eta_k \downarrow 0$. Set $k=0$.
		\item[1.] Apply Algorithm~\ref{alg:opt-fom-str} (in Appendix \ref{app-fom}) to the problem $\min_x \varphi_k(x)$ to  find $x^{k+1}\in\dom(P)$ satisfying
		\begin{equation}
		\label{eq:inext-subprob-m}
		\dist (0, \partial \varphi_k(x^{k+1})) \leq \eta_k,\footnote{In view of Proposition \ref{prop:cplx-fom-str-cvx},  one can terminate Algorithm~\ref{alg:opt-fom-str} once a point $x^{k+1}$ satisfying $2L_{\nabla \phi}\|\tilde{x}^{k+1} - x^{k+1}\| \le \eta_k$ is found, where $\phi(\cdot)= \varphi_k(\cdot)-P(\cdot)$, $\tilde{x}^{k+1}= \prox_{P/L_{\nabla \phi}}\left( x^{k+1} - \nabla \phi(x^{k+1})/L_{\nabla \phi} \right)$, and $L_{\nabla \phi}$ is the Lipschitz constant of $\nabla \phi$. Such $x^{k+1}$ clearly satisfies \eqref{eq:inext-subprob-m}.}
		\end{equation}
		where 
		\begin{equation}
		\label{eq:prox-al-func}
		\varphi_k(x) = \mcL(x,\lambda^k;\rho_k) + \frac{1}{2\rho_k}\|x - x^k\|^2.
		\end{equation}
		\item[2.] Set $\lambda^{k+1} = \Pi_{\mcK^*}\big(\lambda^k + \rho_k g(x^{k+1}) \big)$. 
		\item[3.] If $(x^{k+1},\lambda^{k+1})$ satisfies \eqref{eq:eps-approx-kkt-point} or the following two inequalities are satisfied
		\begin{equation}
		\label{eq:check-stop-m}
		\frac{1}{\rho_k}\|(x^{k+1}, \lambda^{k+1}) - (x^k,\lambda^k)\| \leq \frac{\epsilon}{2}, \quad \eta_k \leq \frac{\epsilon}{2},
		\end{equation}
		output $(x^+, \lambda^{+}) = (x^{k+1}, \lambda^{k+1})$ and terminate the algorithm.
		\item[4.] Set $k \leftarrow k+1$ and go to Step 1.
	\end{itemize}
	{\bf End.}
\end{algo}

For ease of later reference, we refer to the iterations of Algorithm~\ref{alg:opt-fom-str} for solving the 
AL subproblems as the {\it inner iterations} of Algorithm \ref{alg:mi-al}, and call the update from 
$(x^k,\lambda^k)$ to $(x^{k+1},\lambda^{k+1})$ an {\it outer iteration} of Algorithm \ref{alg:mi-al}. 
Notice from \eqref{eq:prox-al-func} that $\varphi_k$ is strongly convex with modulus $1/\rho_k$. The AL subproblem $\min_x \varphi_k(x)$ arising in Algorithm \ref{alg:mi-al} can thus be suitably solved by Algorithm~\ref{alg:opt-fom-str}.

In the rest of this section, we present our main results for Algorithm \ref{alg:mi-al}, whose proof is deferred  to Subsection \ref{sec:pf-mi-al}. Before proceeding, we introduce some further notation that will be used subsequently. Let $X^*$ be the set of optimal solutions of problem \eqref{eq:conic-cp} and $\hat x^*\in X^*$ such that 
$\|x^0-\hat x^*\|=\dist(x^0,X^*)$. In addition, we define 
\beq \label{param-m} 
{\bar D}_X =\max\{D_X,1\}, \ D=  \dist(x^0,X^*)+D_\Lambda, \ \bar{D}= \max\{D,1\}, \ \hat B = L_{\nabla f} + L_{\nabla g}\|\hat\lambda^*\| + L_{\nabla g}D, 
\eeq
where $D_X$, $D_\Lambda$ and $\hat\lambda^*$ are defined in \eqref{param}, and $L_{\nabla f}$ and 
$L_{\nabla g}$ are the Lipschitz constants of $\nabla f$ and $\nabla g$ on $\dom(P)$, respectively. The following theorem shows that Algorithm \ref{alg:mi-al} with a suitable choice of $\{\rho_k\}$ and $\{\eta_k\}$ is guaranteed to find an $\epsilon$-KKT solution of problem  \eqref{eq:conic-cp} in a finite number of outer iterations.

\begin{thm}
	\label{prop:output-gantee-m}
	\bi
	\item[(i)]
	If  Algorithm~\ref{alg:mi-al} successfully terminates, then the output $(x^+,\lambda^+)$  is an $\epsilon$-KKT solution of problem \eqref{eq:conic-cp}.
	\item[(ii)]
	Suppose that $\{\rho_k\}$  and $\{\eta_k\}$ satisfy that
	\beq \label{rho-eta-m}
	\rho_k>0 \ \mbox{is nondecreasing}, \quad 0<\eta_k \downarrow 0, \quad 
	\frac{\sum_{i=0}^{2k}\rho_i\eta_i}{\rho_k\sqrt{k+1}} \to 0.  \footnote{For example,  $\rho_k = \rho_0\alpha^k$ and $\eta_k = \eta_0\beta^k$ satisfy \eqref{rho-eta-m} for any $\rho_0>0$, $\eta_0>0$, $\alpha>1$ and $0<\beta<1/\alpha$.}
	\eeq 
	Then Algorithm~\ref{alg:mi-al} terminates in a finite number of outer iterations. 
	\item[(iii)] Furthermore, if $N$ is a nonnegative integer such that
	\begin{equation}
	\label{eq:cplx-out-general-m}
	\frac{D + \sum_{k=0}^N \rho_k\eta_k}{\rho_N} \leq \frac{\epsilon}{2}, \quad \eta_N \leq \frac{\epsilon}{2},
	\end{equation}
	then Algorithm~\ref{alg:mi-al} terminates in at most $N+1$ outer iterations.
	\ei
\end{thm}       

The next theorem provides an upper bound on the total number of the inner iterations Algorithm~\ref{alg:mi-al}, that is, the total iterations of Algorithm \ref{alg:opt-fom-str} applied to solve all AL subproblems of Algorithm~\ref{alg:mi-al}. 

\begin{thm}
	\label{thm:cplx-mi-al}
	Let $\epsilon>0$ be given, and $\bar D_X$ and $\bar{D}$ be defined in \eqref{param-m}. Suppose that $\{\rho_k\}$ and $\{\eta_k\}$ are chosen as 
	\begin{equation}
	\label{eq:para-mi-al}
	\rho_k = \rho_0\alpha^k, \quad \eta_k = \eta_0\beta^k
	\end{equation}
	for some $\rho_0\geq 1$, $0<\eta_0\leq 1$, $\alpha>1$, $0<\beta<1$ such that $\gamma = \alpha\beta < 1$.
	Then, the total number of inner iterations of Algorithm \ref{alg:mi-al} for finding an $\epsilon$-KKT solution of problem \eqref{eq:conic-cp} is at most
	\begin{equation}
	\label{eq:cplx-mi-al}
	{\cal T}(\epsilon) = \ceil{\frac{8\alpha^2\sqrt{\hat{C}\rho_0 }}{\alpha-1}\log\frac{2\alpha \hat{C}  {\bar D}_X}{\eta_0\beta}}\max\left\{1,\ceil{ \frac{2(\bar{D}+\rho_0\eta_0)}{(1- \gamma)\epsilon}\log_\alpha \frac{2\alpha(\bar{D}+\rho_0\eta_0)}{(1 - \gamma)\epsilon} }\right\},
	\end{equation}
	where $\hat{C}  = C\rho_0 + \hat B + L_{\nabla g}\rho_0\eta_0/(1 - \gamma) + 1$, and $C$ and $\hat B$ are defined in \eqref{bparam} and \eqref{param-m}, respectively.
\end{thm}

\begin{rem} 
One can see from Theorem \ref{thm:cplx-mi-al} that the first-order iteration complexity of Algorithm \ref{alg:mi-al} for finding an $\epsilon$-KKT solution of problem \eqref{eq:conic-cp} is $\cO(\epsilon^{-1}\log \epsilon^{-1})$, which significantly improves the previously best-known iteration-complexity $\mathcal{O}(\epsilon^{-7/4})$ achieved by first-order I-AL methods for finding an $\epsilon$-KKT solution of \eqref{eq:conic-cp}.
\end{rem}

Finally, the following theorem shows that for the $\epsilon$-KKT solution $(x^+,\lambda^+)$ outputted by Algorithm \ref{alg:mi-al}, its primal objective gap, dual objective gap, and constraint violation are at most $\mathcal{O}(\epsilon)$.
\begin{thm}
	\label{thm:eps-funv-m}
	Consider the same setting as in Theorem \ref{thm:cplx-mi-al}. Then, the output $(x^+,\lambda^+)$ of Algorithm \ref{alg:mi-al} satisfies 
\[
	\dist(g(x^+),-\mcK)  \leq \epsilon,  \qquad  -C_1\epsilon \leq F(x^+) - F^*  \leq C_3\epsilon, \qquad 
    0 \leq F^* - d(\lambda^+) \leq (C_1+C_3)\epsilon, 
\]
	where $C_1= \dist(0, \Lambda^*)$ and $C_3 = D_X + \|\hat{\lambda}^*\| + D + \rho_0\eta_0/(1-\gamma)$.
\end{thm}

\section{Proof of the main results}\label{sec:proof}
In this section we prove our main results presented in Sections \ref{sec:i-al} and \ref{sec:mi-al}, that is, Theorems \ref{prop:output-gantee}-\ref{thm:eps-funv-m}. To proceed, we present some technical results that will be used subsequently in our proofs.

Recall that the AL function given in \eqref{eq:al-func} can be written as $\mcL(x,\lambda;\rho) = \mcS(x,\lambda;\rho) + P(x) - \|\lambda^2\|/(2\rho)$, where $\mcS$ is defined in \eqref{S-fun}. The following lemma states some properties of the function $\mcS$, whose proof can be found in Appendix~\ref{proof:al-func-x}.
\begin{lemma}
\label{lemma:al-func-x}
Let $\mcS$ be defined in \eqref{S-fun}. For any $(\lambda,\rho)\in\Re^m\times\Re_{++}$, the following statements hold.
\begin{itemize}
\item[(i)] $\mcS(x,\lambda;\rho)$ is convex and continuously differentiable in $x$ and
\beq \label{grad-S}
 \nabla_x \mcS(x,\lambda;\rho) = \nabla f(x) + \nabla g(x) \Pi_{\mcK^*}\big(\lambda + \rho g(x)\big). 
\eeq
\item[(ii)] $\nabla_x \mcS(x,\lambda;\rho)$ is Lipschitz continuous on $\dom(P)$ with a Lipschitz constant $L$ given by 
$$ L= L_{\nabla f} + L_{\nabla g}\big(\|\lambda\| + \rho M_g\big) + \rho L_g^2. $$
\end{itemize}
\end{lemma}

The lemma below presents some key properties of the AL function, whose proof is a direct  extension of the results in \cite{Rock76b} and thus omitted.

\begin{lemma}
	\label{lemma:al-func-elp}
	For any $(x,\lambda,\rho)\in\Re^n\times\Re^m\times\Re_{++}$, the following identity holds
	\begin{equation*}
	\mcL(x,\lambda;\rho) = \max_{\eta\in\Re^m} \left\{ l(x,\eta) - \frac{1}{2\rho}\|\eta - \lambda\|^2\right\}.
	\end{equation*}
	In addition, if $x\in\dom(P)$, the maximum is attained uniquely at $\bar{\lambda} = \Pi_{\mcK^*}(\lambda + \rho g(x))$. Consequently, the following statements hold.
	\begin{itemize}
		\item[(i)] For any $(\lambda,\rho)\in\Re^m\times\Re_{++}$, $d(\lambda;\rho)$ satisfies that
		\begin{equation}
		\label{eq:al-dual-func-elp-cp}
		d(\lambda;\rho) = \max_{\eta\in\Re^m}\left\{ d(\eta) - \frac{1}{2\rho}\|\eta - \lambda\|^2\right\}. 
		\end{equation}
		\item[(ii)] $\mcL(x,\lambda;\rho)$ is a convex function in $x$, and for any $x\in\dom(P)$, we have
\[
		\partial_x\mcL(x,\lambda;\rho) = \partial_x l(x,\bar{\lambda}).
\]
		\item[(iii)] $\mcL(x,\lambda;\rho)$ is a concave function in $\lambda$, and for any $x\in\dom(P)$, it is differentiable in $\lambda$ and
\[
		\frac{1}{\rho}(\bar{\lambda} - \lambda) = \nabla_\lambda \mcL(x,\lambda;\rho) \in \partial_\lambda l(x,\bar{\lambda}).
\]
	\end{itemize}
\end{lemma}

We next state some results for inexact proximal point algorithm (PPA), whose proof can be found in Appendix \ref{app:proof-ppa}.
\begin{lemma}
\label{lemma:cplx-ippa} 
Let $\mcT:\Re^n\rightrightarrows\Re^n$ be a maximally monotone operator and  $z^*\in\Re^n$  such that $0\in\mcT(z^*)$. Let $\{z^k\}$ be a sequence generated by an inexact PPA, starting with some $z^0$ and obtaining $z^{k+1}$ by approximately evaluating $\mcJ_{\rho_k}(z^k)$ such that 
\begin{equation}
\label{eq:iter-PPA}
\|z^{k+1} - \mcJ_{\rho_k}(z^k)\| \le  e_k
\end{equation}
for some $\rho_k>0$ and $e_k\geq 0$, where $\mcJ_\rho=(\mcI + \rho\mcT)^{-1}$  and $\mcI$ is the identify operator. Then we have
\beqa
\|z^s - z^{*}\| \leq \|z^t - z^{*}\| + \sum_{i=t}^{s-1}e_i, \quad \forall s\geq t\geq 0, \label{eq:iter-ppa-bd} \\
\|z^{k+1} - z^k\| \leq \|z^0 - z^*\| + \sum_{i=0}^k e_i, \quad \forall k\geq 0. \label{eq:ippa-iter-cplx}
\eeqa
Moreover, for any $K\ge 1$, we have
\begin{equation}
\label{eq:ippa-iter-cplx-min}
\min_{K\leq k\leq 2K}\|z^{k+1} - z^k\| \leq \frac{\sqrt{2}\left(\|z^0 - z^{*}\| + 2\sum_{k=0}^{2K}e_k\right)}{\sqrt{K+1}}.
\end{equation}
\end{lemma}

The following two lemmas show that Algorithms~\ref{alg:i-al} and \ref{alg:mi-al}  can be viewed as an inexact PPA method applied to solve the monotone inclusion problems $0\in \mcT_d(\lambda)$ and 
$0\in\mcT_l(x,\lambda)$, respectively.

\begin{lemma}
\label{lemma:i-al=i-ppa}
Let $\{\lambda^k\}$ be the sequence generated by Algorithm~\ref{alg:i-al}. Then for any $k\geq 0$, one has
\[
\|\lambda^{k+1} - \mcJ_{\rho_k}(\lambda^k)\| \leq \sqrt{2\rho_k\eta_k},
\]
where $\mcJ_{\rho_k} = (\mcI + \rho_k\mcT_d)^{-1}$ and $\mcT_d$ is defined in \eqref{eq:def-two-opt1}.
\end{lemma}

\begin{proof}
It follows from the definition of $\dist(\cdot,-\mcK)$ that 
for any $\rho>0$, $\lambda\in\Re^m$ and $x\in\dom(P)$, 
\begin{equation}
\label{eq:dist-reform}
\dist(\lambda + \rho g(x), -\mcK) = \min_u\{\|\lambda - u\|: \rho g(x) + u \preceq_\mcK 0\}.
\end{equation}
By~\cite[Exercise 2.8]{ruszczynski2006nonlinear},  the minimum of \eqref{eq:dist-reform} is attained uniquely at $\bar{u} = \lambda - \Pi_{\mcK^*}\big(\lambda + \rho g(x)\big)$. 
These together with \eqref{eq:al-func} yield
\begin{equation}
\label{eq:alt-def-L}
\mcL(x^{k+1},\lambda^k;\rho_k) = f(x^{k+1}) + P(x^{k+1}) + \frac{1}{2\rho_k}\left[\|\lambda^k - u^k\|^2 - \|\lambda^k\|^2\right], \\
\end{equation}
where $u^k = \lambda^k - \Pi_{\mcK^*}(\lambda^k + \rho_k g(x^{k+1}))$. By this and Step 2 of Algorithm \ref{alg:i-al}, we have $u^k = \lambda^k - \lambda^{k+1}$. Moreover, it follows from \eqref{eq:aug-dual-func} and \eqref{eq:dist-reform} that
\begin{align}
d(\lambda^k;\rho_k) 
= \min_u\left\{ v(u) + \frac{1}{2\rho_k}\left[\|\lambda^k - u\|^2 - \|\lambda^k\|^2\right]\right\}, \label{eq:alt-def-d}
\end{align}
where  
\begin{equation}
\label{eq:def-vu}
v(u) = \min_x \left\{ f(x) + P(x): \rho_k g(x) + u \preceq_\mcK 0 \right\}.
\end{equation}
Since $f+P$ is convex and $g$ is convex with respect to $\mcK$, it is not hard to see that $v$ is also convex. Hence, the objective function in \eqref{eq:alt-def-d} is strongly convex in $u$ and it has a unique minimizer $\bar{u}^k$. Claim that $\bar{u}^k = \lambda^k - \mcJ_{\rho_k}(\lambda^k)$. Indeed, it follows from \eqref{eq:alt-def-d} and Danskin's theorem that $\nabla_\lambda d(\lambda^k;\rho_k) = - \bar{u}^k/\rho_k$. 
In addition, it follows from \eqref{eq:al-dual-func-elp-cp} and the definition of $\mcJ_{\rho_k}(\lambda^k)$ that
$$ d(\lambda^k;\rho_k) = \max_{\eta\in\Re^m} \left\{ d(\eta) - \frac{1}{2\rho_k}\|\eta - \lambda^k\|^2 \right\},$$ 
and the maximum is attained uniquely at $\mcJ_{\rho_k}(\lambda^k)$. By these and Danskin's theorem, we obtain that $\nabla_\lambda d(\lambda^k;\rho_k) = (\mcJ_{\rho_k}(\lambda^k) - \lambda^k)/\rho_k$, which together with $\nabla_\lambda d(\lambda^k;\rho_k) = -\bar{u}^k/\rho_k$ yields $\bar{u}^k = \lambda^k - \mcJ_{\rho_k}(\lambda^k)$ as desired. 
By this, \eqref{eq:alt-def-L}, \eqref{eq:alt-def-d} and \eqref{eq:def-vu}, we obtain that
\begin{align}
\mcL(x^{k+1},\lambda^k;\rho_k) - d(\lambda^k;\rho_k) & = f(x^{k+1}) + P(x^{k+1}) + \frac{1}{2\rho_k}\|\lambda^k - u^k\|^2 - \min_u \left\{ v(u) + \frac{1}{2\rho_k}\|\lambda^k - u\|^2\right\} \nn \\
& \geq v(u^k) + \frac{1}{2\rho_k}\|\lambda^k - u^k\|^2 - \min_u \left\{ v(u) + \frac{1}{2\rho_k}\|\lambda^k - u\|^2\right\} \label{eq:ge-u} \\
& \geq \frac{1}{2\rho_k}\|u^k - \bar{u}^k\|^2 = \frac{1}{2\rho_k}\|\mcJ_{\rho_k}(\lambda^k) - \lambda^{k+1}\|^2, \label{sub-u-bu}
\end{align}
where \eqref{eq:ge-u} follows from \eqref{eq:def-vu} and the fact that 
$$ \rho_k g(x^{k+1}) + u^k = \lambda^k + \rho_k g(x^{k+1}) - \Pi_{\mcK^*}(\lambda^k + \rho_k g(x^{k+1})) = \Pi_{-\mcK}(\lambda^k + \rho_k g(x^{k+1})) \preceq_\mcK 0, $$
and \eqref{sub-u-bu} follows from $\bar{u}^k = \argmin_u\{v(u) + \frac{1}{2\rho_k}\|\lambda^k - u\|^2\}$, the fact that $v(u) + \frac{1}{2\rho_k}\|\lambda^k - u\|^2$ is strongly convex with modulus $1/\rho_k$,  $u^k = \lambda^k - \lambda^{k+1}$, and $\bar{u}^k = \lambda^k - \mcJ_{\rho_k}(\lambda^k)$. The conclusion then follows from \eqref{eq:inext-subprob} and \eqref{sub-u-bu}.
\end{proof}

\begin{lemma}
\label{lemma:mi-al=i-ppa}
Let $\{(x^k,\lambda^k)\}$ be generated by Algorithm~\ref{alg:mi-al}. For any $k\geq 0$, one has
\begin{equation}
\label{eq:mi-al=i-ppa}
\|(x^{k+1},\lambda^{k+1}) - \mcJ_{\rho_k}(x^k,\lambda^k)\| \leq \rho_k\eta_k,
\end{equation}
where $\mcJ_{\rho_k}=(\mcI + \rho_k\mcT_l)^{-1}$ and $\mcT_l$ is defined in \eqref{eq:def-two-opt2}.
\end{lemma}

\begin{proof}
By Lemma~\ref{lemma:al-func-elp} and $\lambda^{k+1} = \Pi_{\mcK^*}\big(\lambda^k + \rho_k g(x^{k+1})\big)$, one has
\begin{equation}
\label{eq:diff-x-lambda-m}
\partial_x \mcL(x^{k+1},\lambda^k;\rho_k) = \partial_x l(x^{k+1},\lambda^{k+1}), \quad \frac{1}{\rho_k}(\lambda^{k+1} - \lambda^k) \in \partial_\lambda l (x^{k+1},\lambda^{k+1}).
\end{equation}
By~\eqref{eq:inext-subprob-m}, there exists $\|v\|\leq \eta_k$ such that
$$ v \in \partial_x \mcL(x^{k+1},\lambda^k;\rho_k) + \frac{1}{\rho_k}(x^{k+1} - x^k). $$
This together with~\eqref{eq:diff-x-lambda-m} implies that
\begin{equation}
\label{eq:diff-x-lambda-m-v}
x^k + \rho_kv \in \rho_k \partial_x l(x^{k+1},\lambda^{k+1}) +  x^{k+1}, \quad \lambda^k \in -\rho_k\partial_\lambda l(x^{k+1}, \lambda^{k+1}) + \lambda^{k+1},
\end{equation}
which, by the definition of $\mcT_l$, are equivalent to 
$$ (x^k + \rho_kv, \lambda^k) \in (\mcI + \rho_k\mcT_l)(x^{k+1},\lambda^{k+1}). $$
It follows from this and $\mcJ_{\rho_k}=(\mcI + \rho_k\mcT_l)^{-1}$ that  $(x^{k+1},\lambda^{k+1}) = \mcJ_{\rho_k}(x^k+\rho_kv,\lambda^k)$. By this and the non-expansion of $\mcJ_{\rho_k}$, we obtain
$$ \|(x^{k+1}, \lambda^{k+1}) - \mcJ_{\rho_k}(x^k,\lambda^k)\| = \|\mcJ_{\rho_k}(x^k + \rho_kv, \lambda^k) - \mcJ_{\rho_k}(x^k,\lambda^k)\| \leq \|\rho_kv\| \leq \rho_k\eta_k, $$
which yields \eqref{eq:mi-al=i-ppa} as desired. 
\end{proof}

\subsection{Proof of the main results for Algorithm \ref{alg:i-al}}\label{sec:pf-i-al}

In this subsection we provide a proof for  Theorems \ref{prop:output-gantee}, \ref{thm:cplx-i-al} and \ref{thm:eps-funv}.

\begin{proof}[Proof of Theorem \ref{prop:output-gantee}] 
(i) One can easily see that $(x^+,\lambda^+)$  is an $\epsilon$-KKT solution of~\eqref{eq:conic-cp} if Algorithm~\ref{alg:i-al} terminates in Step 3.  We now show that it is  also an $\epsilon$-KKT solution of~\eqref{eq:conic-cp} if Algorithm~\ref{alg:i-al} terminates in Step 4. To this end, suppose  that Algorithm~\ref{alg:i-al} terminates in Step 4 at some iteration $k$, that is, the inequalities \eqref{eq:check-stop} hold for some $k$.   For convenience, let $(\tilde{\lambda},\tilde{\rho},\tilde{x})=(\lambda^k,\rho_k,x^{k+1})$. It then follows that $(x^+,\lambda^+) =$ Postprocessing$(\tilde{\lambda},\tilde{\rho},\tilde{x},\epsilon)$,  and \eqref{eq:inext-postpro} and  \eqref{eq:output-postpro} hold for such $\tilde{\lambda}$ and $\tilde{\rho}$. By Definition \ref{defi:approx-KKT}, it suffices to show that $\dist\big(0,\partial_x l(x^+,\lambda^+)\big)\leq \epsilon$ and $\dist\big(0,\partial_\lambda l(x^+,\lambda^+)\big)\leq \epsilon$.

We start by showing $\dist\big(0,\partial_x l (x^+,\lambda^+)\big) \leq \epsilon$. For convenience, let $\varphi(x) = \mcL(x, \tilde{\lambda}; \tilde{\rho})$.  Notice from Lemma~\ref{lemma:al-func-x} that $\nabla_x \mcS(x,\tilde{\lambda};\tilde{\rho})$ is Lipschitz continuous on $\dom(P)$ with Lipschitz constant $\tilde{L}$, where $\tilde L$ is given in \eqref{eq:postpro-para}. Then, by \eqref{eq:postpro-para}, \eqref{eq:inext-postpro}, \eqref{eq:output-postpro} and Proposition~\ref{prop:unconst-cvx-opt} in Appendix \ref{app-fom}, one has $\varphi(x^+) \leq \varphi(\hat{x})$ and
\begin{equation}
\label{eq:primal-eps-1}
\dist\big(0,\partial \varphi(x^+)\big)\leq \sqrt{8\tilde{L}\big(\varphi(\hat{x}) - \min_{x\in\Re^n} \varphi(x)\big)} \leq \sqrt{8\tilde{L}\tilde{\eta}}\leq \epsilon.
\end{equation}
In addition, it follows from \eqref{eq:output-postpro} and  Lemma~\ref{lemma:al-func-elp} that
$$ \partial\varphi(x^+) = \partial_x\mcL(x^+,\tilde{\lambda};\tilde{\rho}) = \partial_x l\big(x^+,\Pi_{\mcK^*}\big(\tilde{\lambda} + \tilde{\rho} g(x^+)\big)\big) = \partial_x l(x^+, \lambda^+). $$
This together with \eqref{eq:primal-eps-1} yields $\dist\big(0,\partial_x l (x^+,\lambda^+)\big) \leq \epsilon$.
It remains to show $\dist\big(0,\partial_\lambda l(x^+,\lambda^+)\big)\leq \epsilon$. Let $\mcJ_{\rho_k} = (\mcI + \rho_k\mcT_d)^{-1}$, where $\mcT_d$ is defined in \eqref{eq:def-two-opt1}. By \eqref{eq:check-stop} 
 and Lemma~\ref{lemma:i-al=i-ppa}, one has 
\[
\|\lambda^{k+1} - \mcJ_{\rho_k}(\lambda^k)\|\leq \sqrt{2\rho_k\eta_k}\leq \frac{\rho_k\epsilon}{8}.
\]
Using this and the first inequality in \eqref{eq:check-stop}, we have 
\[
\|\lambda^k- \mcJ_{\rho_k}(\lambda^k)\|\leq \|\lambda^{k+1}-\lambda^k\| + \|\lambda^{k+1} - \mcJ_{\rho_k}(\lambda^k)\| \le   \frac{3\rho_k\epsilon}{4}+\frac{\rho_k\epsilon}{8} =  \frac{7\rho_k\epsilon}{8},
\]
which, together with $\tilde{\lambda}=\lambda^k$ and $\tilde{\rho}=\rho_k$, leads to $\|\tilde{\lambda}-\mcJ_{\tilde{\rho}}(\tilde{\lambda})\| \le  7\tilde{\rho}\epsilon/8$. In addition, by $\varphi=\mcL(\cdot,\tilde{\lambda};\tilde{\rho})$, the second relation in  \eqref{eq:output-postpro}, and the same arguments as those for \eqref{sub-u-bu}, one has 
$$ \|\lambda^+ - \mcJ_{\tilde{\rho}}(\tilde{\lambda})\| \leq \sqrt{2\tilde{\rho} \big(\mcL(x^+,\tilde{\lambda};\tilde{\rho}) - \min_{x\in\Re^n}\mcL(x,\tilde{\lambda};\tilde{\rho}) \big)} = \sqrt{2\tilde{\rho}\big(\varphi(x^+) - \min_{x\in\Re^n}\varphi(x) \big)}. $$
This, together with $\varphi(x^+) \leq \varphi(\hat{x})$, \eqref{eq:postpro-para} and \eqref{eq:inext-postpro}, yields that
\[
\|\lambda^+ - \mcJ_{\tilde{\rho}}(\tilde{\lambda})\| \leq  \sqrt{2\tilde{\rho}\big(\varphi(\hat{x}) - \min_{x\in\Re^n}\varphi(x) \big)} \leq \sqrt{2\tilde{\rho}\tilde{\eta}} \leq \frac{\tilde{\rho}\epsilon}{8}.
\]
Using this and $\|\tilde{\lambda}-\mcJ_{\tilde{\rho}}(\tilde{\lambda})\| \le  7\tilde{\rho}\epsilon/8$, we obtain that 
\beq \label{eq:dual-resi-bd}
\|\lambda^+ - \tilde{\lambda}\| \le \|\tilde{\lambda} - \mcJ_{\tilde{\rho}}(\tilde{\lambda})\| + \|\lambda^+ - \mcJ_{\tilde{\rho}}(\tilde{\lambda})\|  \le \frac{7\tilde{\rho}\epsilon}{8}+\frac{\tilde{\rho}\epsilon}{8} =\tilde{\rho}\epsilon.
\eeq
Moreover, by Lemma~\ref{lemma:al-func-elp} and the second relation in  \eqref{eq:output-postpro},  one has $(\lambda^+ - \tilde{\lambda})/\tilde{\rho}\in\partial_\lambda l(x^+,\lambda^+)$. 
This along with~\eqref{eq:dual-resi-bd} implies $\dist\big(0,\partial_\lambda l(x^+,\lambda^+)\big)\leq \epsilon$.

(ii) Suppose for contradiction that Algorithm~\ref{alg:i-al} does not terminate in a finite number of iterations. It then follows that \eqref{eq:check-stop} does not hold for any $k$. By Lemma \ref{lemma:i-al=i-ppa} and \eqref{eq:ippa-iter-cplx-min}, one has that for any $k\geq 1$,
\[
\min_{k \leq i\leq 2k} \|\lambda^{i+1} - \lambda^i\| \leq \frac{\sqrt{2}\left(D_\Lambda + 2\sum_{i=0}^{2k}\sqrt{2\rho_i\eta_i}\right)}{\sqrt{k+1}},
\]
where $D_\Lambda$ is defined in \eqref{param}. Since $\{\rho_k\}$ is assumed to be nondecreasing, we further have
\beq \label{lambda-reldiff}
\min_{k \leq i\leq 2k} \frac{1}{\rho_i}\|\lambda^{i+1} - \lambda^i\| \leq \frac{\sqrt{2}\left(D_\Lambda + 2\sum_{i=0}^{2k}\sqrt{2\rho_i\eta_i}\right)}{\rho_k\sqrt{k+1}}.
\eeq
By this and \eqref{rho-eta}, one has that $\min_{k \leq i\leq 2k} \|\lambda^{i+1} - \lambda^i\|/\rho_i \to 0$ and $\eta_k/\rho_k \to 0$ as $k \to \infty$, which imply that \eqref{eq:check-stop} is satisfied for some $k$ and thus leads to a contradiction. 

(iii) Suppose for contradiction that Algorithm~\ref{alg:i-al} does not terminate within $2N+1$ outer iterations. It then follows that \eqref{eq:check-stop} does not hold for all $0\le k \le 2N$.
 Let $\tilde{k}\in\Argmin_{N\leq k\leq 2N} \|\lambda^{k+1} - \lambda^k\|/\rho_k$. By this, \eqref{lambda-reldiff}, and the assumption that $\{\rho_k\}$ is nondecreasing and $\{\eta_k\}$ is decreasing, one has 
$$ \frac{1}{\rho_{\tilde{k}}} \|\lambda^{\tilde{k}+1} - \lambda^{\tilde{k}}\| \leq \frac{\sqrt{2}\left(D_\Lambda + 2\sum_{k=0}^{2N}\sqrt{2\rho_k\eta_k}\right)}{\rho_N\sqrt{N+1}} < \frac{3}{4}\epsilon, \qquad  \frac{\eta_{\tilde{k}}}{\rho_{\tilde{k}}} \leq \frac{\eta_N}{\rho_N}\leq \frac{\epsilon^2}{128}. $$
Hence, \eqref{eq:check-stop} holds for $k = \tilde{k} \leq 2N$, which leads to a contradiction. 
\end{proof}

To prove Theorems \ref{thm:cplx-i-al} and \ref{thm:eps-funv}, we need the following result.

\begin{lemma}
\label{lemma:lip-const-bd}
For any $k\geq 0$, the Lipschitz constant of $\nabla_x\mcS(x,\lambda^k;\rho_k)$, denoted as $L_k$, satisfies
\begin{equation}
\label{eq:lip-const-bd}
L_k \leq C\rho_k + B + L_{\nabla g}\sum_{i=0}^{k-1}\sqrt{2\rho_i\eta_i},
\end{equation}
where $B$ and $C$ are given in \eqref{param} and \eqref{bparam}, respectively.
\end{lemma}
\begin{proof}
By Lemma~\ref{lemma:al-func-x}, one has $L_k \leq L_{\nabla f} + L_{\nabla g}\big(\|\lambda^k\| + \rho_kM_g\big) + \rho_kL_g^2$. Combining \eqref{eq:iter-ppa-bd} with Lemma \ref{lemma:i-al=i-ppa}, and using \eqref{param}, we have
\begin{equation}
\label{eq:lambda-bd}
\|\lambda^k\| \leq \|\hat\lambda^*\| + \|\lambda^k - \hat\lambda^*\| \leq \|\hat\lambda^*\| + D_\Lambda + \sum_{i=0}^{k-1} \sqrt{2\rho_i\eta_i}, 
\end{equation} 
where $\hat\lambda^*$ is defined right above \eqref{param}. 
By these and the definitions of $B$ and $C$, one obtains \eqref{eq:lip-const-bd}. 
\end{proof}

We are now ready to prove Theorems \ref{thm:cplx-i-al} and \ref{thm:eps-funv}.
\begin{proof}[Proof of Theorem \ref{thm:cplx-i-al}]
	For convenience, let $\epsilon_0= \min\{1,\epsilon\}$. Let $\bar{N}$ be the number of outer iterations of Algorithm \ref{alg:i-al}. Also, let $\mcI_k$ and $\mcI_p$ be the number of iterations executed by Algorithm \ref{alg:opt-fom} at the $k$th outer iteration of Algorithm \ref{alg:i-al} and in the subroutine Postprocessing, respectively. In addition, let $T$ be the total number of  inner iterations of Algorithm \ref{alg:i-al}. Clearly, we have $T = \sum_{k=0}^{\bar{N}-1}\mcI_k + \mcI_p. $
	In what follows, we first derive upper bounds on $\bar{N}$, $\mcI_k$ and $\mcI_p$, and then use this formula to obtain an upper bound on $T$.
	
	First, we derive an upper bound on $\bar{N}$. By \eqref{eq:para-i-al}, we have that $\eta_k = \eta_0(k+1)^{-5/2}\sqrt{\epsilon_0}$ for any $k\geq 0$. Hence, for any $K\geq 0$, it holds that
	\begin{equation}
	\label{eq:bd-para}
	\sum_{k=0}^K\sqrt{2\rho_k\eta_k} = \sqrt{2\rho_0\eta_0}\ \epsilon_0^{\frac{1}{4}}\sum_{k=0}^{K} (k+1)^{-\frac{1}{2}}
	\leq 2\sqrt{2\rho_0\eta_0}\ \epsilon_0^{\frac{1}{4}}\sqrt{K+1},
	\end{equation}
	where the inequality follows from $\sum_{k=0}^K(k+1)^{-1/2}\leq 2\sqrt{K+1}$.
	Let	$\gamma= 7\bD^{1/2}\epsilon_0^{-1/2}$ and $N = \ceil{\gamma}$. It follows from~\eqref{eq:para-i-al}, \eqref{eq:bd-para}, and $\gamma \leq N \leq \gamma + 1$ that
	\begin{equation}
	\label{eq:check-outer-bd}
	\frac{D_\Lambda + 2\sum_{k=0}^{2N}\sqrt{2\rho_k\eta_k}}{\rho_{N}\sqrt{N+1}} 
	\leq \frac{\bD + 4\sqrt{2\rho_0\eta_0}\ \epsilon_0^{\frac{1}{4}}\sqrt{2N+1}}{\rho_0 (N+1)^2}
	\leq \frac{\bD}{\rho_0 (N+1)^2} + \frac{8\eta_0^{\frac{1}{2}}\epsilon_0^{\frac{1}{4}}}{\rho_0^{\frac{1}{2}}(N+1)^{\frac32}}. 
	\end{equation}
	Notice that 
	$$ \frac{\bD}{\rho_0 (N+1)^2} \le \frac{\bD}{\rho_0\gamma^2} = \frac{\bD}{\rho_0(49\bD \epsilon_0^{-1})} = \frac{\epsilon_0}{49\rho_0} \leq \frac{\epsilon}{49}, $$
	where the first inequality is by $\gamma \leq N+1$ and the last inequality follows from $\rho_0\geq 1$ and $\epsilon_0\leq \epsilon$.
	Also, by $\bD\geq 1$, we have $\gamma \geq 7\epsilon_0^{-1/2}$. This together with $\gamma\leq N+1$, $\rho_0\geq 1$, and $\eta_0\leq 1$ yields
	$$ 
	\frac{8\eta_0^{\frac{1}{2}}\epsilon_0^{\frac{1}{4}}}{\rho_0^{\frac{1}{2}}(N+1)^{\frac32}}
	\leq \frac{8\epsilon_0^{\frac{1}{4}}}{\gamma^{\frac32}} \leq  \frac{8\epsilon_0^{\frac{1}{4}}}{7^{\frac{3}{2}}\epsilon_0^{-\frac{3}{4}}} = \frac{8\epsilon_0}{7^{\frac32}} < \frac{4\epsilon_0}{9} \leq \frac{4\epsilon}{9},
	$$
	Substituting the above two inequalities into \eqref{eq:check-outer-bd}, one has
	\begin{equation}
	\label{eq:cplx-out-1}
	\frac{D_\Lambda + 2\sum_{k=0}^{2N}\sqrt{2\rho_k\eta_k}}{\rho_{N}\sqrt{N+1}} < \frac{\epsilon}{2}.
	\end{equation}
	In addition, using $N + 1 \geq \gamma \geq  7\epsilon_0^{-1/2}$, \eqref{eq:para-i-al}, $\epsilon_0\leq 1$, $\rho_0\geq 1$ and $\eta_0\leq 1$, we obtain that
	\begin{equation}
	\label{eq:cplx-out-2}
	\frac{\eta_{N}}{\rho_{N}} = \frac{\eta_0\epsilon_0^{\frac{1}{2}}}{\rho_0(N+1)^{4}}\leq \frac{1}{7^4\epsilon_0^{-2}} = \frac{\epsilon_0^2}{7^4} < \frac{\epsilon^2}{128}.
	\end{equation}
	By \eqref{eq:cplx-out-1}, \eqref{eq:cplx-out-2} and Theorem \ref{prop:output-gantee} (iii), we obtain
	\begin{equation}
	\label{eq:N-i-al}
	\bar{N} \leq 2N + 1 = 2\left\lceil 7\bD^{\frac{1}{2}}\epsilon_0^{-\frac{1}{2}}\right\rceil + 1.
	\end{equation}
	
	Second, we derive an upper bound on $\mcI_k$. Let $L_k$ be the Lipschitz constant of $\nabla_x\mcS(x,\lambda^k;\rho_k)$. It follows from \eqref{eq:lip-const-bd}  and \eqref{eq:para-i-al} that for any $k\geq 0$,
	\begin{equation}
	\label{eq:lip-bd-spe}
	L_k \leq \bC\rho_0(k+1)^{\frac{3}{2}} + \bB + 2\sqrt{2\rho_0\eta_0}\ \epsilon_0^{\frac{1}{4}}L_{\nabla g}(k+1)^{\frac{1}{2}}.
	\end{equation}
	This, together with Proposition~\ref{prop:cplx-fom-coro}, ~\eqref{eq:inext-subprob} and~\eqref{eq:para-i-al}, yields that
	\begin{align}
	\mcI_k & \leq \ceil{D_X\sqrt{\frac{2L_k}{\eta_k}}}  \leq 1 + \sqrt{2}D_X\sqrt{\frac{\bC\rho_0(k+1)^{\frac{3}{2}} + \bB + 2\sqrt{2\rho_0\eta_0}\ \epsilon_0^{\frac{1}{4}}L_{\nabla g}(k+1)^{\frac{1}{2}}}{\eta_0(k+1)^{-\frac{5}{2}}\epsilon_0^{\frac{1}{2}}}} \nn \\
	& \leq 1 + \sqrt{2}D_X\sqrt{\frac{\bC\rho_0(k+1)^{\frac{3}{2}} + \bB\rho_0 + 2\sqrt{2}\rho_0\epsilon_0^{\frac{1}{4}}L_{\nabla g}(k+1)^{\frac{1}{2}}}{\eta_0(k+1)^{-\frac{5}{2}}\epsilon_0^{\frac{1}{2}}}} \nn \\
	& \leq 1 + D_X\sqrt{\frac{2\rho_0}{\eta_0}}\left( \bC^{\frac{1}{2}}\epsilon_0^{-\frac{1}{4}}(k+1)^2 + \bB^{\frac{1}{2}}\epsilon_0^{-\frac{1}{4}}(k+1)^{\frac{5}{4}} + 2 L_{\nabla g}^{\frac{1}{2}}\epsilon_0^{-\frac{1}{8}}(k+1)^{\frac{3}{2}} \right), \label{eq:cplx-in-spe}
	\end{align}
	where the third inequality is due to $\rho_0\geq 1$ and $\eta_0\leq 1$, and the last inequality follows by $\sqrt{a+b+c}\leq \sqrt{a} + \sqrt{b} + \sqrt{c}$ for any $a,b,c\geq 0$.
	
	Third, we derive an upper bound on $\mcI_p$. Recall that $\bar{N}$ is the number of outer iterations, that is, \eqref{eq:check-stop} is satisfied when $k = \bar{N}-1$. It then follows from Algorithm \ref{alg:i-al} that $(\tilde{\lambda},\tilde{\rho}) = (\lambda^{\bar{N}-1}, \rho_{\bar{N}-1})$ and $\tilde{L} = L_{\bar{N}-1}$. By these, Proposition \ref{prop:cplx-fom-coro}, \eqref{eq:postpro-para}, \eqref{eq:inext-postpro} and $\epsilon_0\leq \epsilon$, we have
	\begin{equation}
	\label{eq:cplx-pp-spe}
	\mcI_p  \leq \ceil{ D_X\sqrt{ \frac{2L_{\bar{N}-1} }{\tilde\eta} }} \leq \ceil{\frac{16D_X}{\epsilon_0}\cdot\max\left\{\sqrt{\frac{L_{\bar{N}-1}}{\rho_{\bar{N}-1}}}, \frac{L_{\bar{N}-1}}{4}\right\}}
	\end{equation}
	In addition, it follows from \eqref{eq:lip-bd-spe} that
	\begin{equation}
	\label{eq:bd-tildeL}
	L_{\bar{N}-1} \leq \bC\rho_0\bar{N}^{\frac{3}{2}} + \bB + 2\sqrt{2\rho_0\eta_0}\ \epsilon_0^{\frac{1}{4}}L_{\nabla g}\bar{N}^{\frac{1}{2}}.
	\end{equation}
	By this and \eqref{eq:para-i-al}, we obtain that for any  $\bar{N}\geq 1$,
	\begin{align}
	\sqrt{\frac{L_{\bar{N}-1}}{\rho_{\bar{N}-1}}} & \leq \sqrt{\frac{\bC\rho_0\bar{N}^{\frac{3}{2}} + \bB + 2\sqrt{2\rho_0\eta_0}\ \epsilon_0^{\frac{1}{4}}L_{\nabla g}\bar{N}^{\frac{1}{2}}}{\rho_0\bar{N}^{\frac{3}{2}}}} \leq \sqrt{\bC + \bB + 2\sqrt{2}\epsilon_0^{\frac{1}{4}}L_{\nabla g}} \nn \\
	& \leq \bC^{\frac{1}{2}} + \bB^{\frac{1}{2}} + 2\epsilon_0^{\frac{1}{8}}L_{\nabla g}^{\frac{1}{2}}, \label{eq:cplx-pp-spe-1}
	\end{align}
	where the second inequality uses $\bar{N}\geq 1$, $\rho_0\geq 1$ and $\eta_0\leq 1$, and the last inequality follows by $\sqrt{a+b+c}\leq \sqrt{a} + \sqrt{b} + \sqrt{c}$ for any $a,b,c\geq 0$. By \eqref{eq:bd-tildeL}, \eqref{eq:cplx-pp-spe-1}, $\epsilon_0\leq 1$, $\eta_0\leq 1$, $\bC\geq 1$ and $\bB\geq 1$, it is not hard to verify that for all $\bar{N}\geq 1$,
	\begin{equation}
	\label{eq:cplx-pp-spe-2}
	\max\left\{\sqrt{\frac{L_{\bar{N}-1}}{\rho_{\bar{N}-1}}}, \frac{L_{\bar{N}-1}}{4}\right\} \leq 
	\bC\rho_0\bar{N}^{\frac{3}{2}} + \bB + 2\rho_0^{\frac{1}{2}}\epsilon_0^{\frac{1}{8}}\bar{N}^{\frac{1}{2}}(L_{\nabla g}+L_{\nabla g}^{\frac{1}{2}}).
	\end{equation}
	Substituting \eqref{eq:cplx-pp-spe-2} into \eqref{eq:cplx-pp-spe}, we arrive at
	\begin{equation}
	\label{eq:cplx-pp-spe-final}
	\mcI_p \leq  1 + \frac{16D_X}{\epsilon_0}\left( \bC\rho_0\bar{N}^{\frac{3}{2}} + \bB + 2\rho_0^{\frac{1}{2}}\epsilon_0^{\frac{1}{8}}\bar{N}^{\frac{1}{2}}(L_{\nabla g}+L_{\nabla g}^{\frac{1}{2}}).\right).
	\end{equation}
	
	Finally, we use \eqref{eq:N-i-al}, \eqref{eq:cplx-in-spe} and \eqref{eq:cplx-pp-spe-final} to derive an upper bound on the overall complexity $T$. By \eqref{eq:N-i-al}, $N = \ceil{\gamma}$ and $\gamma\geq 7$, one has $\bar{N} - 1 \leq 2N \leq 2\gamma + 2 \leq 3\gamma-1$. This together with \eqref{eq:cplx-in-spe} yields that
	\begin{align}
	\sum_{k=0}^{\bar{N}-1} \mcI_k  & \leq  3\gamma + D_X\sqrt{\frac{2\rho_0}{\eta_0}} \left( \bC^{\frac{1}{2}}\epsilon_0^{-\frac{1}{4}}\sum_{k=0}^{\floor{3\gamma}-1} (k+1)^2 + \bB^{\frac{1}{2}}\epsilon_0^{-\frac{1}{4}}\sum_{k=0}^{\floor{3\gamma}-1} (k+1)^{\frac{5}{4}} + 2L_{\nabla g}^{\frac{1}{2}}\epsilon_0^{-\frac{1}{8}}\sum_{k=0}^{\floor{3\gamma}-1} (k+1)^{\frac{3}{2}} \right) \nn \\
	& \leq 3\gamma + D_X\sqrt{\frac{2\rho_0}{\eta_0}} \left( \frac{8}{3}\bC^{\frac{1}{2}}\epsilon_0^{-\frac{1}{4}}(3\gamma)^3 + \frac{2^{\frac{17}{4}}}{9}\bB^{\frac{1}{2}}\epsilon_0^{-\frac{1}{4}}(3\gamma)^{\frac{9}{4}} + \frac{2^{\frac92}}{5}L_{\nabla g}^{\frac{1}{2}}\epsilon_0^{-\frac{1}{8}} (3\gamma)^{\frac{5}{2}} \right)  \nn \\
	& \leq 3\gamma + 72D_X\sqrt{\frac{2\rho_0}{\eta_0}}\left(\bC^{\frac{1}{2}}\epsilon_0^{-\frac{1}{4}}\gamma^3 + \bB^{\frac{1}{2}}\epsilon_0^{-\frac{1}{4}}\gamma^{\frac{9}{4}} + L_{\nabla g}^{\frac{1}{2}}\epsilon_0^{-\frac{1}{8}}\gamma^{\frac{5}{2}} \right), \nn
	\end{align}
	where the second inequality is due to 
	\[
	\sum_{k=0}^{K-1} (k+1)^\alpha \leq \frac{1}{1+\alpha}(K+1)^{1 + \alpha} \le \frac{2^{1+\alpha}}{1+\alpha}K^{1 + \alpha}, \quad \forall \alpha>0, K\geq 1.
	\]
	Recall that $\gamma = 7\bD^{1/2}\epsilon_0^{-1/2}$. Substituting this into the above inequality, we obtain
	\[
	\sum_{k=0}^{\bar{N}-1} \mcI_k = \mathcal{O}\left( \frac{D_X\bD^{3/2}\bC^{1/2}}{\epsilon_0^{7/4}} + \frac{D_X\bD^{9/8}\bB^{\frac{1}{2}} + D_X\bD^{5/4}L_{\nabla g}^{1/2}}{\epsilon_0^{11/8}} + \frac{\bD^{1/2}}{\epsilon_0^{1/2}}\right).
	\]
	In addition, by $\bar{N}\leq 3\gamma$, $\gamma = 7\bD^{1/2}\epsilon_0^{-1/2}$ and \eqref{eq:cplx-pp-spe-final}, we obtain that
	\[
	\mcI_p = \mathcal{O}\left( \frac{D_X\bD^{3/4}\bC}{\epsilon_0^{7/4}} + \frac{D_X\bB}{\epsilon_0} + \frac{D_X\bD^{1/4}(L_{\nabla g}+L_{\nabla g}^{1/2})}{\epsilon_0^{9/8}}\right).
	\]
	Recall that $T = \sum_{k=0}^{\bar{N}-1}\mcI_k + \mcI_p$. By these, $\bD\geq 1$, $\bC\geq 1$, and $\bB\geq 1$, we have
	$$ T = \mathcal{O}\left( \frac{D_X\bD^{3/2}\bC}{\epsilon_0^{7/4}} + \frac{D_X\bD^{5/4}\bB^{1/2}(1+L_{\nabla g}^{1/2})}{\epsilon_0^{11/8}} + 
	\frac{D_X\bD^{1/4}(L_{\nabla g} + L_{\nabla g}^{1/2})}{\epsilon_0^{9/8}} + \frac{D_X\bB}{\epsilon_0} + \frac{\bD^{1/2}}{\epsilon_0^{1/2}} \right). $$
	This together with $\epsilon_0 = \min\{1,\epsilon\}$ yields the complexity bound in Theorem \ref{thm:cplx-i-al}.
\end{proof}

\begin{proof}[Proof of Theorem \ref{thm:eps-funv}]
Recall from Theorem \ref{prop:output-gantee} that the output $(x^+,\lambda^+)$ of Algorithm \ref{alg:i-al} is an $\epsilon$-KKT solution of \eqref{eq:conic-cp}. It then follows from Remark \ref{rem1} that $\dist(g(x^+), \mcN_{\mcK^*}(\lambda^+))\leq \epsilon$. It also follows from Definition \ref{defi:approx-KKT} that there exist ${u}\in\partial_x l(x^+,\lambda^+)$ and ${v}\in\partial_\lambda l(x^+,\lambda^+)$ with $\|{u}\|\leq \epsilon$ and $\|{v}\|\leq \epsilon$. 

We start by proving the first inequality in \eqref{eq:primal-val-bd}. It is not hard to verify that $\mcN_{\mcK^*}(\lambda^+)\subseteq -\mcK$, which together with  $\dist(g(x^+), \mcN_{\mcK^*}(\lambda^+))\leq \epsilon$ implies that $\dist(g(x^+),-\mcK) \leq \epsilon$ holds as desired. 

We next prove the second inequality in \eqref{eq:primal-val-bd}. Let $\lambda^*\in\Lambda^*$ be such that $\|{\lambda}^*\|=\dist(0,\Lambda^*)$. Also, let $x^*$ be any optimal solution of problem \eqref{eq:conic-cp}. Then, $(x^*,{\lambda}^*)$ is a saddle point of $l$, which implies that $l(x^+,{\lambda}^*) - l(x^*,{\lambda}^*) \geq 0$ and $F^* = l(x^*,{\lambda}^*)$.
By these and $l(x^+,{\lambda}^*) = F(x^+) + \langle {\lambda}^*, g(x^+)\rangle$, we obtain
\begin{equation}
\label{eq:bd-3}
F(x^+) - F^* \geq -\langle {\lambda}^*, g(x^+) \rangle. 
\end{equation} 
In addition, by ${\lambda}^*\in\mcK^*$ and $\dist(g(x^+),-\mcK) \leq \epsilon$, one has that
$$ \langle {\lambda}^*, g(x^+) \rangle = \langle {\lambda}^*, g(x^+) - \Pi_{-\mcK}(g(x^+))\rangle + \underbrace{\langle {\lambda}^*,\Pi_{-\mcK}(g(x^+))\rangle}_{\leq 0} \leq \|{\lambda}^*\|\cdot\dist(g(x^+),-\mcK) \leq \|{\lambda}^*\|\epsilon,$$
which, along with  \eqref{eq:bd-3} and $C_1=\dist(0,\Lambda^*)=\|{\lambda}^*\|$, implies that $F(x^+) - F^* \geq -C_1\|\epsilon$.

We now prove the last inequality in \eqref{eq:primal-val-bd}. Recall that $u\in\partial_x l(x^+,\lambda^+)$ and $v\in\partial_\lambda l(x^+,\lambda^+)$ with $\|u\|\leq \epsilon$ and $\|v\|\leq \epsilon$. By \eqref{eq:lag-dual-func}, $\lambda^+\in\mcK^*$ and Assumption \ref{ass:function}, there exists an $\tilde{x}\in\dom(P)$ satisfying $d(\lambda^+) = l(\tilde{x},\lambda^+)$. These, together with the convexity of $l(\cdot,\lambda^+)$, $x^+\in\dom(P)$ and \eqref{param}, imply that 
\begin{equation}
\label{eq:bd-1}
l(x^+,\lambda^+) - d(\lambda^+) = l(x^+,\lambda^+) - l(\tilde{x},\lambda^+) \leq \langle x^+ - \tilde{x}, u \rangle \leq \|x^+ - \tilde{x}\|\|u\| \leq D_X\epsilon.
\end{equation}
Also, using ${v}\in\partial_\lambda l(x^+,\lambda^+)$ and \eqref{eq:express-Tl}, we obtain that $g(x^+) - {v}\in\mcN_{\mcK^*}(\lambda^+)$. Then, by the definition of $\mcN_{\mcK^*}(\lambda^+)$ and the fact that $\mcK^*$ is a closed convex cone, it is not hard to see that 
$\langle g(x^+) - v, \lambda^+ \rangle = 0$. Thus, we obtain
\begin{equation}
\label{eq:bd-11}
\langle \lambda^+, g(x^+)\rangle = \langle \lambda^+, v \rangle \geq -\|\lambda^+\|\|v\| \geq -\epsilon\|\lambda^+\|.
\end{equation} 
Let $\bar{N}$ be the number of outer iterations of Algorithm \ref{alg:i-al}.
Claim that 
\beq \label{eq:bd-2}
\|\lambda^+\| \leq \rho_{\bar{N}-1}\epsilon + \|\hat{\lambda}^*\| + D_\Lambda + \sum_{i=0}^{\bar{N}-1}\sqrt{2\rho_i\eta_i}.
\eeq
Indeed, due to \eqref{eq:lambda-bd}, \eqref{eq:bd-2} clearly holds if Algorithm \ref{alg:i-al} terminates in Step 3. We now assume that it terminates in Step 4. Then $(x^+,\lambda^+) =$ Postprocessing$(\lambda^{\bar{N}-1},\rho_{\bar{N}-1},x^{\bar{N}},\epsilon)$. 
 This together with \eqref{eq:dual-resi-bd} and \eqref{eq:lambda-bd} yields that
\[
\|\lambda^+\| \leq \|\lambda^+ -\lambda^{\bar{N}-1}\| + \|\lambda^{\bar{N}-1}\| \leq \rho_{\bar{N}-1}\epsilon + \|\hat{\lambda}^*\| + D_\Lambda + \sum_{i=0}^{\bar{N}-2}\sqrt{2\rho_i\eta_i},
\]
and hence \eqref{eq:bd-2} holds.  By \eqref{eq:para-i-al}, \eqref{eq:N-i-al}, $\bar{D}_\Lambda\geq 1$, and $\epsilon_0 = \min\{1,\epsilon\}$, we obtain
$$ \rho_{\bar{N}-1} \epsilon = \rho_0\bar{N}^{\frac{3}{2}}\epsilon \leq \rho_0\left(17\bar{D}_\Lambda^{\frac{1}{2}}\epsilon_0^{-\frac{1}{2}} \right)^{\frac{3}{2}}\epsilon \leq 71\rho_0\bar{D}_\Lambda\cdot\max\{1,\epsilon\}. $$
Also, by \eqref{eq:bd-para}, \eqref{eq:N-i-al}, $\bar{D}_\Lambda\geq 1$, $\rho_0\geq 1$ and $\eta_0\leq 1$, we obtain
$$ \sum_{i=0}^{\bar{N}-1}\sqrt{2\rho_i\eta_i} \leq 2\sqrt{2\rho_0\eta_0}\epsilon_0^{\frac{1}{4}}\sqrt{\bar{N}} \leq 2\sqrt{2\rho_0\eta_0}\epsilon_0^{\frac{1}{4}}\sqrt{17\bar{D}_\Lambda^{\frac{1}{2}}\epsilon_0^{-\frac{1}{2}}} \leq 12\sqrt{\rho_0\eta_0}\bar{D}_\Lambda^{\frac{1}{4}} \leq 12\rho_0\bar{D}_\Lambda. $$
These, together with \eqref{eq:bd-2}, $\rho_0\ge 1$ and $D_\Lambda\leq \bar{D}_\Lambda$,  yield that
$
\|\lambda^+\| \leq C_2\cdot\max\{1,\epsilon\},
$
where $C_2 = 84\rho_0\bar{D}_\Lambda + \|\hat{\lambda}^*\|$. By this, \eqref{eq:bd-1}, \eqref{eq:bd-11}, and $l(x^+,\lambda^+) = F(x^+) + \langle\lambda^+, g(x^+)\rangle$, we obtain
\begin{equation}
\label{eq:primal-dual-fun-gap}
F(x^+) - d(\lambda^+) \leq D_X\epsilon - \langle \lambda^+, g(x^+)\rangle \leq (D_X + \|\lambda^+\|)\epsilon \leq (D_X+ C_2\cdot\max\{1,\epsilon\})\epsilon,
\end{equation} 
which along with $d(\lambda^+)\leq F^*$ results in the last inequality in \eqref{eq:primal-val-bd}. 

Finally, the first inequality in \eqref{eq:dual-val-bd} trivially holds, while the second inequality in \eqref{eq:dual-val-bd} follows from \eqref{eq:primal-dual-fun-gap} and the second inequality in \eqref{eq:primal-val-bd}. 
\end{proof}

\subsection{Proof of the main results for Algorithm \ref{alg:mi-al}}\label{sec:pf-mi-al}

In this subsection we provide a proof for Theorems \ref{prop:output-gantee-m}, \ref{thm:cplx-mi-al}, and \ref{thm:eps-funv-m}. 

\begin{proof}[Proof of Theorem \ref{prop:output-gantee-m}]
(i) Suppose that Algorithm~\ref{alg:mi-al} terminates in Step 3 at some iteration $k$. It then follows 
that $(x^+,\lambda^+)$  is already an $\epsilon$-KKT solution of problem \eqref{eq:conic-cp} 
or  the inequalities \eqref{eq:check-stop-m} hold for such $k$.  We next show that for the latter case, $(x^+, \lambda^+) = (x^{k+1},\lambda^{k+1})$ is also an $\epsilon$-KKT solution of~\eqref{eq:conic-cp}. To this end, suppose that \eqref{eq:check-stop-m} holds at the $k$th iteration. Notice that \eqref{eq:diff-x-lambda-m-v} holds for some $\|v\|\leq \eta_k$, which yields  
\[
\frac{1}{\rho_k}(x^k-x^{k+1}) + v \in \partial_x l(x^{k+1},\lambda^{k+1}), \quad \frac{1}{\rho_k}(\lambda^{k+1} - \lambda^k) \in \partial_\lambda l (x^{k+1},\lambda^{k+1}).
\]
By these, \eqref{eq:check-stop-m} and $(x^+, \lambda^+) = (x^{k+1},\lambda^{k+1})$, we obtain
\begin{align*}
\mbox{dist}(0,\partial_x l(x^+, \lambda^+)) & \leq \frac{1}{\rho_k}\|x^{k+1} - x^k - \rho_kv\| \leq \frac{1}{\rho_k}\|x^{k+1} - x^k\| + \|v\| \leq \epsilon,\\
\mbox{dist}(0,\partial_\lambda l(x^+, \lambda^+)) & \leq \frac{1}{\rho_k}\|\lambda^{k+1} - \lambda^k\| \leq \epsilon,
\end{align*}
which along with Definition \ref{defi:approx-KKT} imply that $(x^+,\lambda^+)$ is an $\epsilon$-KKT solution of problem \eqref{eq:conic-cp}.

(ii) Suppose for contradiction that Algorithm~\ref{alg:mi-al} does not terminate in a finite number of iterations. It then follows that  \eqref{eq:check-stop-m} does not hold for any $k$. By \eqref{eq:ippa-iter-cplx-min} and Lemma \ref{lemma:mi-al=i-ppa}, one has that
$$ \min_{k \leq i\leq 2k} \|(x^{i+1},\lambda^{i+1}) - (x^i,\lambda^i)\| \leq \frac{\sqrt{2}\left(D + 2\sum_{i=0}^{2k}\rho_i\eta_i\right)}{\sqrt{k+1}}, $$ 
where $D$ is defined in \eqref{param-m}. Since $\{\rho_k\}$ is assumed to be nondecreasing, we further have
\[
\min_{k \leq i\leq 2k} \frac{1}{\rho_i}\|(x^{i+1},\lambda^{i+1}) - (x^i,\lambda^i)\| \leq \frac{\sqrt{2}\left(D + 2\sum_{i=0}^{2k}\rho_i\eta_i\right)}{\rho_k\sqrt{k+1}}.
\]
By this and \eqref{rho-eta-m}, one has that $\min_{k \leq i\leq 2k} \|(x^{i+1},\lambda^{i+1}) - (x^i,\lambda^i)\|/\rho_i \to 0$ and $\eta_k\to 0$ as $k \to \infty$, which implies that \eqref{eq:check-stop-m} is satisfied for some $k$ and thus leads to a contradiction. 

(iii) Let $\hat x^*$ and $\hat\lambda^*$ be defined right above \eqref{param} and \eqref{param-m}, respectively. Suppose for contradiction that Algorithm~\ref{alg:mi-al} does not terminate within $N+1$ outer iterations. It then follows that \eqref{eq:check-stop-m} does not hold for all $0\le k \le N$. Then, by \eqref{param-m}, \eqref{eq:ippa-iter-cplx}, and \eqref{eq:mi-al=i-ppa}, we have that for all $0\le k \le N$,
\begin{equation*}\label{eq:mi-ppa-bd}
\|(x^{k+1},\lambda^{k+1}) - (x^k,\lambda^k)\| \leq \|(x^0,\lambda^0) - (\hat x^*, \hat\lambda^*)\| + \sum_{i=0}^k\rho_i\eta_i \leq D + \sum_{i=0}^k\rho_i\eta_i,
\end{equation*}
where $D$ is given in \eqref{param-m}. 
By this and \eqref{eq:cplx-out-general-m}, one can see that \eqref{eq:check-stop-m} is satisfied when $k = N$, which leads to a contradiction.  
\end{proof}

To prove Theorems \ref{thm:cplx-mi-al} and \ref{thm:eps-funv-m}, we need the following result.

\begin{lemma}
\label{lem:lip-const-bd-m}
Let $s_k(x) = \mcS(x,\lambda^k;\rho_k) + \|x - x^k\|^2/(2\rho_k)$. Then $s_k$ is continuously differentiable, and moreover, $\nabla s_k$ is Lipschitz continuous on $\dom(P)$ with a Lipschitz constant $L_k$ given by
\begin{equation}
\label{eq:lip-const-bd-m}
L_k = C\rho_k + \hat B + L_{\nabla g}\sum_{i=0}^{k-1} \rho_i\eta_i + \rho^{-1}_k,
\end{equation}
where $C$ and $\hat B$ are defined in \eqref{bparam} and \eqref{param-m}, respectively.
\end{lemma}
\begin{proof}
By the definition of $s_k(x)$ and Lemma~\ref{lemma:al-func-x}, one has
\beq \label{lipschtiz-sk}
\|\nabla s_k(x)-\nabla s_k(y)\| \leq \left(L_{\nabla f} + L_{\nabla g}(\|\lambda^k\| + \rho_kM_g) + \rho_kL_g^2 + \rho^{-1}_k\right) \|x-y\|, \quad \forall x,y \in \dom(P). 
\eeq
In addition, it follows from Lemmas \ref{lemma:cplx-ippa} and \ref{lemma:mi-al=i-ppa} that 
\[
\|(x^k,\lambda^k)- (\hat x^*, \hat\lambda^*)\|\leq \|(x^0,\lambda^0) - (\hat x^*, \hat\lambda^*)\| + \sum_{i=0}^{k-1} \rho_i\eta_i,
\]
where $\hat x^*$ and $\hat\lambda^*$ are defined right above \eqref{param} and \eqref{param-m}, respectively. Hence, we have that
\begin{align}
\|\lambda^k\| \leq \|\hat\lambda^*\| + \|\lambda^k - \hat\lambda^*\| \leq \|\hat\lambda^*\| + D + \sum_{i=0}^{k-1}\rho_i\eta_i, \label{eq:lambda-bd-m}
\end{align}
where $D$ is given in \eqref{param-m}. Substituting this into \eqref{lipschtiz-sk}, and using 
the definitions of $\hat B$ and $C$, we obtain that $\|\nabla s_k(x)-\nabla s_k(y)\| \leq L_k \|x-y\|$ for all  
$x,y \in \dom(P)$. Hence, the conclusion holds.
\end{proof}

We are now ready to provide a proof for Theorems \ref{thm:cplx-mi-al} and \ref{thm:eps-funv-m}.

\begin{proof}[Proof of Theorem \ref{thm:cplx-mi-al}]
Let $\bar{N}$ be the number of outer iterations of Algorithm \ref{alg:mi-al}, and let $\mcI_k$ be the number of first-order iterations executed by  Algorithm \ref{alg:opt-fom-str} at the $k$th outer iteration of Algorithm \ref{alg:mi-al}. In addition, let $T$ be the total number of first-order inner iterations of Algorithm \ref{alg:mi-al}. Clearly, we have $T = \sum_{k=0}^{\bar{N}-1}\mcI_k$.
In what follows, we first derive upper bounds on $\bar{N}$ and $\mcI_k$,  and then use this formula to obtain an upper bound on $T$.

We first derive an upper bound on $\bar{N}$.  Due to~\eqref{eq:para-mi-al} and $0<\gamma<1$, we have that 
\begin{equation}
\label{eq:bd-para-m}
\sum_{k=0}^K\rho_k\eta_k = \rho_0\eta_0\sum_{k=0}^K\gamma^k \leq \rho_0\eta_0\sum_{k=0}^\infty\gamma^k = \frac{\rho_0\eta_0}{1 - \gamma}, \quad \forall K \ge 0.
\end{equation}
Let 
$$ N = \max\left\{1,  \ceil{\log_\alpha\frac{2(\bar{D}+\rho_0\eta_0)}{(1-\gamma)\epsilon}}\right\}. $$
Since $N \geq \log_\alpha\frac{2(\bar{D}+\rho_0\eta_0)}{(1 - \gamma)\epsilon}$, we have from \eqref{eq:para-mi-al} that
$$ \rho_N \geq \frac{2\rho_0(\bar{D}+\rho_0\eta_0)}{(1-\gamma)\epsilon}. $$
By this, \eqref{eq:bd-para-m}, $D\leq\bar{D}$, and $\rho_0\geq 1$, we obtain
$$ \frac{D + \sum_{k=0}^N\rho_k\eta_k}{\rho_N} \leq \frac{\bar{D} + \frac{\rho_0\eta_0}{1 - \gamma}}{\frac{2\rho_0(\bar{D}+\rho_0\eta_0)}{(1 - \gamma)\epsilon}} = \frac{\epsilon}{2}\cdot\frac{\bar{D}(1 - \gamma)+\rho_0\eta_0}{\rho_0(\bar{D} + \rho_0\eta_0)} \leq \frac{\epsilon}{2}\cdot\frac{\bar{D} + \rho_0\eta_0}{\bar{D} + \rho_0\eta_0} = \frac{\epsilon}{2}. $$
In addition, one can observe that $1< \alpha <\beta^{-1}$ and $\bar{D}+\rho_0\eta_0 \ge 1 - \gamma$. By these,  we have 
$$N \geq \log_\alpha\frac{2(\bar{D}+\rho_0\eta_0)}{(1 - \gamma)\epsilon} \geq \log_{\beta^{-1}}\frac{2}{\epsilon},$$ which, together with \eqref{eq:para-mi-al}, $\beta<1$ and $\eta_0\leq 1$, implies that $\eta_N \leq \epsilon/2$. It then follows from these and Theorem \ref{prop:output-gantee-m} (iii) that
\begin{equation}
\label{eq:cplx-out-spe-m}
\bar{N} \leq N + 1 \leq \max\left\{1, \ceil{\log_\alpha\frac{2(\bar{D} + \rho_0\eta_0)}{(1-\gamma)\epsilon}}\right\} + 1.
\end{equation}

We next derive an upper bound on $\mcI_k$. By~\eqref{eq:para-mi-al}, \eqref{eq:lip-const-bd-m},  $\alpha > 1$ and $\rho_0\geq 1$, one has that for any $k\geq 0$,
$$ L_k \leq C\rho_0\alpha^k + \hat B + \frac{L_{\nabla g}\rho_0\eta_0}{1 - \gamma} + \frac{1}{\rho_0\alpha^k} \leq \hat{C} \alpha^k, $$
where $\hat{C}  = C\rho_0 + \hat B + L_{\nabla g}\rho_0\eta_0/(1 - \gamma) + 1$. Notice that $\varphi_k(x)$ is strongly convex with modulus $\mu_k=1/\rho_k$. By this, \eqref{eq:inext-subprob-m}, $\rho_k = \rho_0\alpha^{k}$, $\hat{C} \geq 1$, ${\bar D}_X\geq 1$, $\alpha>1$, $\beta<1$, $\rho_0\geq 1$, $\eta_0\leq 1$, and Proposition~\ref{prop:cplx-fom-str-cvx} (in Appendix \ref{app-fom}), we obtain that for any $k\geq 0$, 
\begin{align}
\mcI_k &\leq \left\lceil{\sqrt{\frac{L_k}{\mu_k}}}\ \right\rceil  \max \left\{ 1, \ceil{2\log \frac{2L_kD_X}{\eta_k}} \right\} \leq \ceil{\sqrt{\hat{C}\rho_0 }\ \alpha^k}\max\left\{ 1, \ceil{ 2\log \frac{2\alpha^k\hat{C}  {\bar D}_X}{\eta_0\beta^k}} \right\} \nn \\
 & \leq \ceil{\sqrt{\hat{C}\rho_0 }\ \alpha^k} \ceil{ 2\log \frac{2\alpha^k\hat{C} {\bar D}_X}{\eta_0\beta^k}}  \leq \left( \sqrt{\hat{C}\rho_0 }\ \alpha^k + 1\right) \left( 2\log \frac{2\alpha^k\hat{C} {\bar D}_X}{\eta_0\beta^k} + 1\right) \nn \\
& \leq 8\sqrt{\hat{C}\rho_0 }\ \alpha^k \log \frac{2\alpha^k\hat{C} {\bar D}_X}{\eta_0\beta^k} \leq 8\sqrt{\hat{C}\rho_0 }\ k\alpha^k \log\frac{2\alpha \hat{C}  {\bar D}_X}{\eta_0\beta}, \label{eq:cplx-in-spe-m} 
\end{align}
 where the third and fifth inequalities follow from $\sqrt{\hat{C}\rho_0 }\ \alpha^k \geq 1$ and $2\log \frac{2\alpha^k \hat{C}  {\bar D}_X}{\eta_0\beta^k} \geq 2\log 2 \geq 1$.

Finally, we derive an upper bound on $T$. By \eqref{eq:cplx-in-spe-m}, one has
\begin{equation}
\label{eq:series-bd-m} 
T = \sum_{k=0}^{\bar{N}-1} \mcI_k \leq 8\sqrt{\hat{C}\rho_0 }\log\frac{2\alpha \hat{C}  {\bar D}_X}{\eta_0\beta} \sum_{k=0}^{\bar{N}-1}k\alpha^k \leq \frac{8\sqrt{\hat{C}\rho_0 }}{\alpha-1}\log\frac{2\alpha \hat{C}  {\bar D}_X}{\eta_0\beta} (\bar{N}-1)\alpha^{\bar{N}},
\end{equation}
where the last inequality is due to $\sum_{k=0}^K k\alpha^k \leq K\alpha^{K+1}/(\alpha - 1)$ for any $K\geq 0$. We divide the rest of the proof into the following two cases.  

Case (a): $\frac{2(\bar{D}+\rho_0\eta_0)}{(1 - \gamma)\epsilon} \geq \alpha$. This along with \eqref{eq:cplx-out-spe-m} implies that $\bar{N} \leq \log_\alpha \frac{2(\bar{D}+\rho_0\eta_0)}{(1 - \gamma)\epsilon} +2$.  By this and  \eqref{eq:series-bd-m}, one has 
$$ T \leq \frac{8\sqrt{\hat{C}\rho_0 }}{\alpha-1}\log\frac{2\alpha \hat{C} {\bar D}_X}{\eta_0\beta} \log_\alpha \frac{2\alpha(\bar{D}+\rho_0\eta_0)}{(1 - \gamma)\epsilon} \cdot \frac{2\alpha^2(\bar{D}+\rho_0\eta_0)}{(1- \gamma)\epsilon}. $$

Case (b): $\frac{2(\bar{D}+\rho_0\eta_0)}{(1 - \gamma)\epsilon} < \alpha$. This together with \eqref{eq:cplx-out-spe-m}  implies that $\bar{N} \leq 2$. By this and \eqref{eq:series-bd-m}, one has 
$$ T \leq \frac{8\alpha^2\sqrt{\hat{C}\rho_0 }}{\alpha-1}\log\frac{2\alpha \hat{C}  {\bar D}_X}{\eta_0\beta}. $$

Combining the results in the above two cases, we obtain \eqref{eq:cplx-mi-al} as desired.
\end{proof}

\begin{proof}[Proof of Theorem \ref{thm:eps-funv-m}]
Let $\bar{N}$ be the number of outer iterations of Algorithm \ref{alg:mi-al}. 
By \eqref{eq:lambda-bd-m} and $\lambda^+ = \lambda^{\bar{N}}$, one has $\|\lambda^+\| \leq \|\hat{\lambda}^*\| + D + \sum_{i=0}^{\bar{N}-1}\rho_i\eta_i$, where $\hat\lambda^*$ is defined right above \eqref{param-m} and $D$ is given in \eqref{param-m}.
This together with \eqref{eq:para-mi-al} and \eqref{eq:bd-para-m} yields that $\|\lambda^+\| \leq \|\hat{\lambda}^*\| + D + \rho_0\eta_0/(1 - \gamma)$.
The rest of the proof follows the same arguments as those in the proof of Theorem \ref{thm:eps-funv}.
\end{proof}

\section{Numerical results}\label{sec:numerical}

In this section we conduct some preliminary numerical experiments to test the performance of our proposed algorithms (Algorithms \ref{alg:i-al} and \ref{alg:mi-al}), and compare them with a closely related I-AL method and its modified version proposed in \cite{Lan16}, which are named as I-AL$_1$ and I-AL$_2$ respectively for ease of reference. In particular, we apply all these algorithms to the linear programming (LP) problem
\begin{equation}
\label{eq:lp-exp}
\min_{x\in\Re^n} \left\{ c^Tx: Ax = b, \ l\leq x\leq u\right\}
\end{equation}
for some $A\in\Re^{m\times n}$, $b\in\Re^m$, $c\in\Re^n$, and $l,u\in\Re$. It is clear that \eqref{eq:lp-exp} is a special case of problem \eqref{eq:conic-cp} with $f(x) = c^Tx$, $P$ being the indicator function of the set $\{x\in\Re^n: l\leq x\leq u\}$, $g(x) = Ax - b$, and $\mcK = \{0\}^m$. All the algorithms are coded in Matlab and all the computations are performed on a Dell desktop with a 3.40-GHz Intel Core i7-3770 processor and 16~GB of RAM. 

In our experiment, we choose $\epsilon = 0.01$ for all the aforementioned algorithms. In addition, the parameters $\{\rho_k\}$ and $\{\eta_k\}$ of these algorithms are set as follows. For Algorithm \ref{alg:i-al}, we set them by \eqref{eq:para-i-al} with $\rho_0 = 100$ and $\eta_0 = 1$. For Algorithm \ref{alg:mi-al}, we choose them by \eqref{eq:para-mi-al} with $\rho_0 = 100$, $\eta_0 = 0.1$, $\alpha=1.1$ 
and $\beta=0.8$. For the algorithms I-AL$_1$ and I-AL$_2$, we choose $\{\rho_k\}$ and $\{\eta_k\}$ as described in \cite{Lan16} and set $t_0 = 1$ as the initial value in their ``guess-and-check'' procedures. 

We randomly generate 20 instances for problem \eqref{eq:lp-exp}, each of which is generated by a similar manner as described in \cite{LLM11}. In particular, given positive integers $m<n$ and a scalar $0<\zeta\leq 1$, we first randomly generate a matrix $A\in\Re^{m\times n}$ with density $\zeta$, whose entries are randomly chosen from the standard normal distribution.\footnote{The matrix $A$ is generated via the Matlab command \textsf{A = sprandn(m,n,$\zeta$)}.} We then generate a vector $x\in\Re^n$ with entries randomly chosen from the uniform distribution on $[-5,5]$ and set $b = Ax$. Also, we generate a vector $c\in\Re^n$ with entries randomly chosen from the standard normal distribution. Finally, we randomly choose $l$ and $u$ from the uniform distribution on $[-10,-5]$ and $[5,10]$, respectively.

The computational results of all the algorithms for solving problem \eqref{eq:lp-exp} with the above
 20 instances are presented in Table \ref{tab:result}. In detail, the parameters $n$, $m$, and $\zeta$ of each instance are listed in the first three columns, respectively. For each instance, the total number of first-order iterations and the CPU time (in seconds) for these algorithms are given in the next four columns and the last four columns, respectively. One can observe that Algorithm \ref{alg:mi-al} performs best in terms of both number of iterations and CPU time, which is not surprising as it has the lowest first-order iteration complexity $\cO(\epsilon^{-1}\log\epsilon^{-1})$ among these algorithms. In addition, although Algorithm~\ref{alg:i-al} and I-AL$_1$ share the same order of first-order iteration complexity $\cO(\epsilon^{-7/4})$, one can observe that the practical performance of Algorithm~\ref{alg:i-al} is substantially better than that of I-AL$_1$. The main reason is perhaps that Algorithm \ref{alg:i-al} uses the dynamic $\{\rho_k\}$ and $\{\eta_k\}$, while I-AL$_1$ uses the static ones through all iterations and also needs a ``guess-and-check'' procedure for approximating the unknown parameter $D_\Lambda$. Finally, we observe that I-AL$_2$ performs much better than I-AL$_1$ and generally better than Algorithm \ref{alg:i-al}, but it is substantially outperformed by Algorithm \ref{alg:mi-al}.

\renewcommand{\arraystretch}{1.2}
\begin{table}[t!]
	\centering
	\small
	\caption{Computational results for solving problem \eqref{eq:lp-exp}}
	\smallskip
	
	\label{tab:result}
	\makebox[\textwidth]{
		\begin{tabular}{|c c c||c c c c || c c c c|}
			\hline
			\multicolumn{3}{|c||}{Parameters} & \multicolumn{4}{|c||}{Iterations ($\times 10^3$)} & 
			\multicolumn{4}{|c|}{CPU Time (in seconds)} \\
			$n$ & $m$ & $\zeta$ & Algorithm \ref{alg:i-al}  & Algorithm \ref{alg:mi-al} & I-AL$_1$ & I-AL$_2$ & Algorithm \ref{alg:i-al}  & Algorithm \ref{alg:mi-al} & I-AL$_1$ & I-AL$_2$ \\
			\hline
			1,000 & 100 & 0.01 & 5 & 13 & 164 & 52 & 0.7 & 0.9 & 18.8 & 6.6 \\
			1,000 & 100 & 0.05 & 8 & 13 & 200 & 23 & 1.2 & 1.2 & 31.5 & 3.8 \\
			1,000 & 100 & 0.10 & 8 & 16 & 200 & 25 & 1.8 & 2.0 & 41.7 & 5.4 \\
			1,000 & 500 & 0.01 & 22 & 16 & 200 & 30 & 3.8 & 1.7 & 33.7 & 5.3 \\
			1,000 & 500 & 0.05 & 23 & 19 & 300 & 35 & 10.8 & 6.3 & 136.9 & 16.5 \\
			1,000 & 500 & 0.10 & 22 & 15 & 300 & 22 & 17.5 & 8.9 & 237.2 & 17.0 \\
			1,000 & 900 & 0.01 & 150 & 20 & 900 & 77 & 35.2 & 3.0 & 208.0 & 18.6 \\
			1,000 & 900 & 0.05 & 124 & 19 & 1,100 & 64 & 94.3 & 10.7 & 876.0 & 51.8 \\
			1,000 & 900 & 0.10 & 132 & 21 & 600 & 49 & 197.2 & 23.9 & 903.3 & 71.0 \\
			5,000 & 500 & 0.01 & 19 & 27 & 200 & 78 & 17.2 & 13.6 & 181.0 & 74.0 \\
			5,000 & 500 & 0.05 & 20 & 31 & 200 & 49 & 46.5 & 49.9 & 505.1 & 126.9 \\
			5,000 & 500 & 0.10 & 19 & 26 & 200 & 42 & 129.9 & 149.6 & 1,357.3 & 288.3 \\
			5,000 & 2,500 & 0.01 & 79 & 20 & 300 & 49 & 225.8 & 40.5 & 852.1 & 140.7 \\
			5,000 & 2,500 & 0.05 & 80 & 27 & 300 & 61 & 1,706.4 & 505.1 & 6,406.2 & 1,309.8 \\
			5,000 & 2,500 & 0.10 & 81 & 31 & 300 & 54 & 3,577.7 & 1,240.9 & 13,324.2 & 2,530.2 \\
			5,000 & 4,500 & 0.01 & 400 & 27 & 1,400 & 191 & 2,953.1 & 167.9 & 10,364.8 & 1,425.8 \\
			5,000 & 4,500 & 0.05 & 406 & 29 & 1,300 & 207 & 17,724.6 & 1,067.8 & 55,608.2 & 8,812.9 \\
			5,000 & 4,500 & 0.10 & 300 & 32 & 1,200 & 172 & 26,489.9 & 2,449.3 & 104,523.0 & 15,002.9 \\
			10,000 & 1,000 & 0.01 & 27 & 30 & 200 & 54 & 76.7 & 52.2 & 572.8 & 157.0 \\
			10,000 & 5,000 & 0.01 & 116 & 29 & 400 & 111 & 1,988.5 & 406.6 & 6,895.0 & 1,931.0 \\
			\hline	
		\end{tabular}
	}
\end{table}

\section{Concluding remarks}\label{sec:conclude}
In this paper our analyses of the I-AL methods rely on the assumption that the domain of the function $P$ is 
compact. One natural question is whether this assumption can be dropped.  In addition, can 
the first-order iteration complexity $\mathcal{O}(\epsilon^{-1}\log\epsilon^{-1})$ for finding an $\epsilon$-KKT solution of problem \eqref{eq:conic-cp} be further improved for 
an I-AL method? These will be left for the future research. 
 
\appendix

\section{An example about the dependence of $D^\epsilon_\Lambda$ on $\epsilon$}\label{app-example}

In this part we present an example to demonstrate $D^\epsilon_\Lambda=\Theta(1/\epsilon)$, where $D^\epsilon_\Lambda=\min\{\|\lambda^0-\lambda\|:\lambda\in\Lambda^*_\epsilon\}$ and $\Lambda^*_\epsilon$ is the set of optimal Lagrangian multipliers associated with the perturbed problem \eqref{eq:conic-lp-pert}.

Consider the linear programming problem:
\begin{equation}
\label{eq:LP1}
\begin{aligned}
& \min_{x\in X} \quad x_1 + x_2 - \delta x_3 \\
& \,\ \mbox{s.t.} \quad x_1 = 1,  \ \delta x_2 = -\delta
\end{aligned}
\end{equation}
for some $\delta\in(0,1)$, where  
$$ X = \left\{x\in\Re^3: -x_2 + \frac{\delta}{1-\delta}x_3 \leq \frac{1}{1 - \delta}, \ -2\leq x_1,x_2,x_3\leq 2 \right\}. $$
We also consider a perturbed problem for \eqref{eq:LP1} given by
\begin{equation}
\label{eq:LP1-p}
\begin{aligned}
& \min_{x\in X} \quad x_1 + x_2 - \delta x_3 + \frac{\epsilon}{4D_X}\|x - x^0\|^2\\
& \,\ \mbox{s.t.} \quad x_1 = 1, \ \delta x_2 = -\delta,
\end{aligned}
\end{equation}
where $x^0 = (1,-1,-1)^T \in X$, $\epsilon = 2D_X\delta$, and $D_X = \max\{\|x - y\|: x,y\in X\}$. It is clear that \eqref{eq:LP1} and \eqref{eq:LP1-p} are a special instance of problems \eqref{eq:conic-lp} and \eqref{eq:conic-lp-pert}, respectively. In addition, one can verify that for any $\delta\in(0,1)$, $\bar{x}^\epsilon = (1,-1,0)^T$ is the optimal solution of \eqref{eq:LP1-p} and $\bar{\lambda}^\epsilon = (1,2D_X/\epsilon)^T$ is the unique optimal Lagrangian multiplier associated with the constraints $x_1 = 1$ and $\delta x_2 = -\delta$ of \eqref{eq:LP1-p}. Assume without loss of generality that $\lambda^0 = 0$. Then, for any $\delta\in(0,1)$, we have 
$D_\Lambda^\epsilon = \|\lambda^0 - \bar{\lambda}^\epsilon\| = \epsilon^{-1}\sqrt{\epsilon^2 + 4D_X^2}$, 
whose dependence on $\epsilon$ is roughly $\Theta(\epsilon^{-1})$ when $\epsilon$ is small. 

\section{Optimal first-order methods for unconstrained convex optimization problems}\label{app-fom}
In this part we review optimal first-order methods for solving a class of convex optimization problems in the form of
\begin{equation}\label{eq:uncons-cvx-opt}
\Psi^* = \min_{x\in\Re^n} \left\{\Psi(x) := \phi(x) + h(x) \right\}, 
\end{equation}
where $\phi,h:\Re^n\rightarrow(-\infty,\infty]$ are closed convex functions, $\phi$ is continuously differentiable on an open set containing $\dom(h)$, and $\nabla \phi$ is Lipschitz continuous with Lipschitz constant $L_{\nabla \phi}$ on $\dom(h)$. In addition, assume that 
$\dom(h)$ is compact and let $D_h := \max_{x,y\in\dom(h)}\|x-y\|$.

We first state a property of problem \eqref{eq:uncons-cvx-opt}, which is used in the proof of some main results of this paper. Its proof follows from some standard arguments and is omitted due to page limit.
\begin{prop}
	\label{prop:unconst-cvx-opt}
	For any $x\in\dom(h)$, we have $\Psi(x^+) \leq \Psi(x)$ and
	\begin{equation*}
	\label{eq:subdiff-func-value}
	\dist\big(0,\partial \Psi(x^+)\big) \leq 2L_{\nabla \phi}\|x^+ - x\| \le \sqrt{8L_{\nabla \phi}\big(\Psi(x) - \Psi^*\big)},
	\end{equation*}
	where $x^+= \prox_{h/L_{\nabla \phi}}\big(x - \nabla \phi(x)/L_{\nabla \phi}\big)$. 
\end{prop}

We next present an optimal first-order method for solving \eqref{eq:uncons-cvx-opt} in which $\phi$ is convex but not necessarily strongly convex. It is a variant of Nesterov's optimal first-order methods \cite{Nesterov05,Nesterov13} and has been studied in, for example, \cite[Section 3]{Tseng08}.

\begin{algo}[An optimal first-order method for \eqref{eq:uncons-cvx-opt} with general convex $\phi$]
	\label{alg:opt-fom}
	\normalfont
	\mbox{}
	\begin{itemize}
		\item[0.] Input $x^0 = z^0 \in \mbox{dom}(h)$. Set $k=0$.
		\item[1.] Set $y^k = (k x^k + 2z^k)/(k+2)$.
		\item[2.] Compute $z^{k+1}$ as
		\[
		z^{k+1} = \argmin_{z}\Big\{\ell(z;y^k) + \frac{L_{\nabla \phi}}{k+2}\|z - z^k\|^2\Big\},
		\]
		where
		\begin{equation}
		\label{eq:def-ell-phi}
		\ell(x;y):= \phi(y) + \Inner{\nabla \phi(y)}{x - y} + h(x).
		\end{equation}
		\item[3.] Set $x^{k+1} = (k x^k + 2z^{k+1})/(k+2)$.  
		\item[4.] Set $k\leftarrow k+1$ and go to Step 1.
	\end{itemize}
	{\bf End.}
\end{algo}

The following result provides an iteration-complexity of Algorithm \ref{alg:opt-fom} for finding an $\epsilon$-optimal solution of \eqref{eq:uncons-cvx-opt}. It is an immediate consequence of \cite[Corollary 1]{Tseng08} and its proof is thus omitted.
\begin{prop}
	\label{prop:cplx-fom-coro}
	Let $\{(x^k,y^k)\}$ be generated by Algorithm \ref{alg:opt-fom} and $\ell(\cdot;\cdot)$ be defined in \eqref{eq:def-ell-phi}. Then, $\Psi(x^k) - \Psi^*\leq \Psi(x^k)-\underline\Psi_k$ for all $k \ge 1$. Moreover, for any given $\epsilon>0$, Algorithm \ref{alg:opt-fom} finds an approximate solution $x^k$ of problem \eqref{eq:uncons-cvx-opt} such that $\Psi(x^k) - \Psi^*\leq \Psi(x^k) - \underline\Psi_k \le \epsilon$ in no more than $K(\epsilon)$ iterations, where
	\begin{equation}\label{eq:cplx-bd-fom}
	K(\epsilon) = \ceil{D_h\sqrt{\frac{2L_{\nabla \phi}}{\epsilon}}},  \qquad \underline\Psi_k = \frac{4}{k(k+2)}\min_x \left\{\sum_{i=0}^{k-1}\frac{i+2}{2}\ell(x;y^i)\right\}.
	\end{equation}
\end{prop}

\begin{rem} \label{phik}
Observe from \eqref{eq:cplx-bd-fom} that 
$
\underline\Psi_k =\frac{2}{k(k+2)}\left[v_k +\min_x \left\{\langle u^k,x\rangle+h(x)\right\}\right],
$
where
\[
u^k=\sum_{i=0}^{k-1}(i+2)\nabla \phi(y^i), \quad v_k = \sum_{i=0}^{k-1}(i+2) \left(\phi(y^i)- \langle \nabla \phi(y^i), y^i\rangle\right).
\]
Note that $(u^k,v_k)$ can be recursively and thus cheaply computed. Once $(u^k,v_k)$ is available, computing $\underline\Psi_k$ only requires solving the problem $\min_x \left\{\langle u^k,x\rangle+h(x)\right\}$, which typically has a closed-form solution.
\end{rem}

We now turn to consider the case of problem \eqref{eq:uncons-cvx-opt} in which $\phi$ is strongly convex, that is, there exists a constant $\mu\in(0,L_{\nabla \phi})$ such that
\begin{equation}
\label{eq:str-cvx}
\Inner{\nabla \phi(x) - \nabla \phi(y)}{x - y} \geq \mu\|x - y\|^2, \quad \forall x,y\in\dom(h).
\end{equation}
We next propose a slight variant of Nesterov's optimal method \cite{Nesterov13,LLX15} for solving problem \eqref{eq:uncons-cvx-opt} with a strongly convex $\phi$.

\begin{algo}[An optimal first-order method for \eqref{eq:uncons-cvx-opt} with strongly convex $\phi$]
	\label{alg:opt-fom-str}
	\normalfont
	\mbox{}
	\begin{itemize}
		\item[0.] Input $x^{-1}\in\dom(h)$, $L_{\nabla \phi}>0$ and $0<\mu< L_{\nabla \phi}$. Compute 
		\begin{equation}
		\label{eq:pgm-extra}
		x^0 = \prox_{h/L_{\nabla \phi}} \left(x^{-1} - \frac{1}{L_{\nabla \phi}}\nabla \phi(x^{-1})\right).
		\end{equation}
		Set $z^0 = x^0$, $\alpha = \sqrt{\mu/L_{\nabla \phi}}$ and $k=0$.
		\item[1.] Set $y^k = (x^k + \alpha z^k)/(1 + \alpha)$.
		\item[2.] Compute $z^{k+1}$ as
		\begin{equation*}
		z^{k+1} = \argmin_{z}\left\{\ell(z;y^k) + \frac{\alpha L_{\nabla \phi}}{2}\|z - \alpha y^k - (1 - \alpha)z^k\|^2\right\},
		\end{equation*}
		where $\ell(x;y)$ is defined in \eqref{eq:def-ell-phi}.
		\item[3.] Set $x^{k+1} = (1 - \alpha) x^k + \alpha z^{k+1}$.
		\item[4.] Set $k\leftarrow k+1$ and go to Step 1.
	\end{itemize}
	{\bf End.}
\end{algo}

\begin{rem} 
	Algorithm \ref{alg:opt-fom-str} differs from Nesterov's optimal method \cite{Nesterov13,LLX15} in that it executes a proximal step \eqref{eq:pgm-extra} to generate $x^0$ while the latter method simply 
	sets $x^0=x^{-1}$.
\end{rem}

We next state an iteration-complexity result for Algorithm \ref{alg:opt-fom-str} for finding an approximate optimal solution of problem \eqref{eq:uncons-cvx-opt}. Its proof follows from  \cite[Theorem 1]{LLX15}, Proposition \ref{prop:unconst-cvx-opt}  and some standard arguments, and is omitted due to page limit.
\begin{prop}
	\label{prop:cplx-fom-str-cvx}
	Suppose that \eqref{eq:str-cvx} holds. Let $\{x^k\}$ be the sequence generated by Algorithm~\ref{alg:opt-fom-str} and $\tilde{x}^k= \prox_{h/L_{\nabla \phi}}\left( x^{k} - \nabla \phi(x^{k})/L_{\nabla \phi} \right)$ 
	for all $k \ge 0$. Then, $\dist(0, \partial \Psi(\tilde{x}^k)) \leq 2L_{\nabla \phi}\|\tilde{x}^k - x^k\|$ for all $k \ge 0$. Moreover, for any given $\epsilon>0$, an approximate solution $\tilde{x}^k$ of problem \eqref{eq:uncons-cvx-opt} satisfying $\dist(0, \partial \Psi(\tilde{x}^k)) \leq 2L_{\nabla \phi}\|\tilde{x}^k - x^k\| \le \epsilon$ is generated by running Algorithm~\ref{alg:opt-fom-str} for at most $\widetilde K(\epsilon)$ iterations, where
	\begin{equation*}\label{eq:tilde-K}
	\widetilde K(\epsilon) = \left\lceil \sqrt{\frac{L_{\nabla \phi}}{\mu}}\ \right\rceil \max\left\{ 1, \left\lceil 2\log \frac{2L_{\nabla \phi} D_h}{\epsilon}\right\rceil\right\}.
	\end{equation*}
\end{prop}

\section{Proof of Lemma \ref{lemma:al-func-x}}\label{proof:al-func-x}
\begin{proof}
	(i) We first show that $\mcS(\cdot,\lambda;\rho)$ is convex. Let $x,x^\prime\in\Re^n$ and $\alpha\in[0,1]$ be arbitrarily chosen. Using \eqref{eq:conic-cvx} and the relation
	\[
	\lambda + \rho g\big(\alpha x + (1-\alpha) x^\prime\big) = \lambda + \rho[\alpha g(x) + (1 - \alpha) g(x^\prime)] + \rho \left(g\big(\alpha x + (1-\alpha) x^\prime\big) - [\alpha g(x) + (1 - \alpha) g(x^\prime)]\right),
	\]
we see that $\lambda + \rho\big[\alpha g(x) + (1 - \alpha) g(x^\prime)\big] \preceq_{-\mcK}  \lambda + \rho g\big(\alpha x + (1-\alpha) x^\prime\big)$.
By this and the fact that $\mcK$ is a convex cone, it is not hard to show that
	\beq \label{dist-ineq}
	\dist^2\Big(\lambda + \rho g\big(\alpha x + (1-\alpha) x^\prime\big), -\mcK\Big) \leq \dist^2\Big(\lambda + \rho[\alpha g(x) + (1 - \alpha) g(x^\prime)], -\mcK\Big).
	\eeq
	In addition, by the convexity of $\dist^2(\cdot,-\mcK)$, one has
	\beqas
	\dist^2\Big(\lambda + \rho[\alpha g(x) + (1 - \alpha) g(x^\prime)], -\mcK\Big) 
	&=&  \dist^2\Big(\alpha \big(\lambda + \rho g(x)\big) + (1 - \alpha) \big(\lambda+\rho g(x^\prime)\big), -\mcK\Big)  \\
	&\leq&  \alpha\cdot \dist^2\big(\lambda + \rho g(x), -\mcK\big) + (1 - \alpha)  \dist^2\big(\lambda + \rho g(x^\prime), -\mcK\big),
	\eeqas
	which along with \eqref{dist-ineq} leads to 
	\[
	\dist^2\Big(\lambda + \rho g\big(\alpha x + (1-\alpha) x^\prime\big), -\mcK\Big) \leq   \alpha 
\cdot\dist^2\big(\lambda + \rho g(x), -\mcK\big) + (1 - \alpha)  \dist^2\big(\lambda + \rho g(x^\prime), -\mcK\big).
	\]
	It thus follows that $\dist^2(\lambda + \rho g(\cdot), -\mcK) $ is convex. This together with the convexity of $f$ implies that $\mcS(\cdot,\lambda;\rho)$ is convex. Next we show that $\mcS(\cdot,\lambda;\rho)$ is continuously differentiable. By the definition of $\dist(\cdot,-\mcK)$, one has 
	$$ \mcS(x,\lambda;\rho) = f(x) + \frac{1}{2\rho} \min_{v\in-\mcK}\|\lambda + \rho g(x) - v\|^2, $$
	where the minimum is attained uniquely at $v = \Pi_{-\mcK}\big(\lambda + \rho g(x) \big)$. 
	Using Danskin's theorem (e.g., see~\cite{Bertsekas99}), we conclude that $\mcS(x,\lambda;\rho)$ is differentiable in $x$ and 
	$$ \nabla_x \mcS(x,\lambda;\rho) = \nabla f(x) + \nabla g(x) \big[ \lambda + \rho g(x) - \Pi_{-\mcK}\big(\lambda + \rho g(x) \big) \big] = \nabla f(x) + \nabla g(x) \Pi_{\mcK^*}\big(\lambda + \rho g(x)\big), $$
where the second equality follows from~\cite[Exercise 2.8]{ruszczynski2006nonlinear}.
	
	(ii) Recall that $\nabla f$, $\nabla g$ and $g$ are Lipschitz continuous on $\dom(P)$. By this and 
	\eqref{grad-S}, we have that for any $x,x^\prime\in\dom(P)$, 
	\begin{align*}
	& \|\nabla_x \mcS(x, \lambda ;\rho) - \nabla_x \mcS(x^\prime, \lambda; \rho)\| = 
	\|\nabla f(x) + \nabla g(x) \Pi_{\mcK^*}\big(\lambda + \rho g(x)\big)-\nabla f(x') - \nabla g(x') \Pi_{\mcK^*}\big(\lambda + \rho g(x')\big)\| \\
	& \ \leq \|\nabla g(x)\Pi_{\mcK^*}\big(\lambda + \rho g(x)\big) - \nabla g(x^\prime)\Pi_{\mcK^*}\big(\lambda + \rho g(x^\prime)\big)\| + \|\nabla f(x) - \nabla f(x^\prime)\| \\
	& \ \leq L_{\nabla g}\|x - x^\prime\|\|\Pi_{\mcK^*}\big(\lambda + \rho g(x)\big)\| + \|\nabla g(x^\prime)\|\|\Pi_{\mcK^*}\big(\lambda + \rho g(x)\big) - \Pi_{\mcK^*}\big(\lambda + \rho g(x^\prime)\big)\| + L_{\nabla f}\|x - x^\prime\|\\ 
	& \ \leq L_{\nabla g}\|x - x^\prime\|\|\lambda + \rho g(x)\| + \rho L_g\|g(x) - g(x^\prime)\| + L_{\nabla f}\|x - x^\prime\| \\
	& \ \leq \big(L_{\nabla g}(\|\lambda\| + \rho M_g) + \rho L_g^2 + L_{\nabla f}\big)\|x - x^\prime\|
	\end{align*}
	where the third inequality is due to the non-expansiveness of the projection operator $\Pi_{\mcK^*}$ and $\|\nabla g(x^\prime)\|\leq L_g$, and the last one follows from  $\|g(x)\|\leq M_g$ and the Lipschitz continuity of $g$ on $\dom(P)$.
\end{proof}

\section{Proof of Lemma \ref{lemma:cplx-ippa}}\label{app:proof-ppa}
\begin{proof}
Since $0\in\mcT(z^*)$ and $\mcT$ is maximally monotone, it follows from \cite[Proposition 1]{Rock76a} that
\begin{equation}
\label{eq:suff-dec-ppa-1}
\| \mcJ_\rho(z) - z^*\|^2 + \|\mcJ_\rho(z) - z\|^2 \leq \|z - z^*\|^2,
\end{equation}
which implies that 
\begin{equation}
\label{eq:suff-dec-ppa-2}
\|\mcJ_\rho(z) - z^*\| \leq \|z - z^*\|, \quad \|\mcJ_\rho(z) - z\| \leq \|z - z^*\|, \quad \forall z \in \Re^n.
\end{equation}
	Let $\xi^k=z^{k+1} - \mcJ_{\rho_k}(z^k)$ for all $k\ge 0$. By this, \eqref{eq:suff-dec-ppa-2}, and \eqref{eq:iter-PPA} with $\rho = \rho_k$ and $z = z^k$, one has
	\[
	\|z^{k+1} - z^{*}\| \leq \|z^{k+1} - \mcJ_{\rho_k}(z^k)\| + \|\mcJ_{\rho_k}(z^k) - z^{*}\| 
	\leq \|\xi^k\| + \|z^k - z^{*}\|, \quad \forall k \ge 0.
	\]
	Summing up the above inequality from $k=t$ to $k=s-1$ yields 
	\beq \label{bdd-1}
	\|z^s - z^{*}\| \leq \|z^t - z^{*}\| + \sum_{i=t}^{s-1} \|\xi^i\|, \quad \forall s\geq t\geq 0.
	\eeq
	Notice from \eqref{eq:iter-PPA} that $\|\xi^k\| \le e_k$ for all $k\ge 0$, which along with \eqref{bdd-1} leads to \eqref{eq:iter-ppa-bd}. Besides, by~\eqref{eq:suff-dec-ppa-2} with $\rho=\rho_k$ and $z = z^k$, one has $ \|\mcJ_{\rho_k}(z^k) - z^k\|  \leq \|z^k - z^*\|$. These together with \eqref{eq:iter-PPA} and \eqref{eq:iter-ppa-bd} yield
	$$ \|z^{k+1} - z^k\| \leq \|z^{k+1} - \mcJ_{\rho_k}(z^k)\| + \|\mcJ_{\rho_k}(z^k) - z^k\| \leq e_k + \|z^k - z^*\| \leq \|z^0 - z^*\| + \sum_{i=0}^k e_i.  $$
	In addition, by the definition of $\xi^k$, and \eqref{eq:suff-dec-ppa-1} 
	with $\mcJ = \mcJ_{\rho_k}$ and $z = z^k$, one has
	\begin{equation*}
	\begin{split}
	\|\mcJ_{\rho_k}(z^k) - z^k\|^2 & \leq \|z^k - z^{*}\|^2 - \|\mcJ_{\rho_k}(z^k) - z^{*}\|^2 
	= \|z^k - z^{*}\|^2 - \|\mcJ_{\rho_k}(z^k) - z^{k+1} + z^{k+1} - z^{*}\|^2 \\
	& \leq \|z^k - z^{*}\|^2 - \|z^{k+1} - z^{*}\|^2 - \|\xi^k\|^2 + 2\|\xi^k\|\|z^{k+1} - z^{*}\|.
	\end{split}
	\end{equation*}
	Summing up the above inequality from $k = K$ to $k = 2K$ and using \eqref{bdd-1}, we obtain that 
	\begin{align}
	\sum_{k=K}^{2K}\|\mcJ_{\rho_k}(z^k) - z^k\|^2 & \leq \|z^K - z^{*}\|^2 - \sum_{k=K}^{2K}\|\xi^k\|^2 + 2\sum_{k=K}^{2K}\|\xi^k\|\left(\|z^K - z^{*}\| + \sum_{j=K}^k\|\xi^j\|\right) \nn\\
	& = \|z^K - z^{*}\|^2 - \sum_{k=K}^{2K}\|\xi^k\|^2 + 2\|z^K - z^{*}\|\cdot\sum_{k=K}^{2K}\|\xi^k\| + 2\sum_{k=K}^{2K}\sum_{j=K}^k\|\xi^k\|\|\xi^j\| \nn\\
	& = \|z^K - z^{*}\|^2 - \sum_{k=K}^{2K}\|\xi^k\|^2 + 2\|z^K - z^{*}\|\cdot\sum_{k=K}^{2K}\|\xi^k\| + \sum_{k=K}^{2K}\|\xi^k\|^2 + \left(\sum_{k=K}^{2K}\|\xi^k\|\right)^2 \nn\\
	& = \|z^K - z^{*}\|^2 + 2\|z^K - z^{*}\|\cdot\sum_{k=K}^{2K}\|\xi^k\| + \left(\sum_{k=K}^{2K}\|\xi^k\|\right)^2 \nn\\
	& = \left( \|z^K - z^{*}\| + \sum_{k=K}^{2K}\|\xi^k\| \right)^2 \leq \left( \|z^0 - z^{*}\| + \sum_{k=0}^{2K}\|\xi^k\| \right)^2, \label{eq:bound-sum-jz-z}
	\end{align}
	where \eqref{eq:bound-sum-jz-z} follows from \eqref{bdd-1} with $t=0$ and $s = K$. Again, by the definition of $\xi^k$, one has
	$$ \|z^{k+1} - z^k\|^2 = \|\mcJ_{\rho_k}(z^k) + \xi^k- z^k \|^2 \leq 2\left(\|\mcJ_{\rho_k}(z^k) - z^k\|^2 + \|\xi^k\|^2\right). $$
	This together with~\eqref{eq:bound-sum-jz-z} yields
	\begin{align*}
	\sum_{k = K}^{2K}\|z^{k+1} - z^k\|^2 & \leq 2\sum_{k = K}^{2K}\|\mcJ_{\rho_k}(z^k) - z^k\|^2 + 2\sum_{k = K}^{2K}\|\xi^k\|^2 
	 \leq 2\left( \|z^0 - z^{*}\| + \sum_{k=0}^{2K}\|\xi^k\| \right)^2 + 2\sum_{k=0}^{2K}\|\xi^k\|^2 \\
	& \leq 2\left( \|z^0 - z^{*}\| + 2\sum_{k=0}^{2K}\|\xi^k\| \right)^2,
	\end{align*}
	which along with $\|\xi^k\| \le e_k$ leads to \eqref{eq:ippa-iter-cplx-min}. The proof is then completed.
\end{proof}

\end{document}